\documentclass[draftcls,onecolumn]{IEEEtran}          

\usepackage{url,amsmath,amsfonts,amssymb,epsfig,epstopdf, color,accents,psfrag}
\usepackage[tight,footnotesize]{subfigure}
\usepackage{lineno}
\usepackage{anysize}

\usepackage{enumerate}

\usepackage[colorinlistoftodos, textwidth=2.5cm, shadow]{todonotes}

\title{ Optimal parameter selection for the alternating direction method of multipliers (ADMM): quadratic problems}

\author{Euhanna Ghadimi,  Andr\'{e} Teixeira, Iman Shames, and Mikael Johansson 
\thanks{
E.~Ghadimi, A.~Teixeira, and M.~Johansson are with the ACCESS Linnaeus Center, Electrical Engineering, Royal Institute of Technology, Stockholm, Sweden.
{\tt\small \{euhanna, andretei, mikaelj\}@ee.kth.se}.
I.~Shames is with the Department of Electrical and Electronic Engineering, The University of Melbourne, Melbourne, Australia.
{\tt\small iman.shames@unimelb.edu.au}. This work was sponsored in part by the Swedish Foundation for Strategic Research, SSF, and the Swedish Research Council, VR.
}
}

\begin{document}
\newtheorem{theorem}{Theorem}
\newtheorem{proposition}{Proposition}
\newtheorem{lem}{Lemma}
\newtheorem{rem}{Remark}
\newtheorem{assumption}{Assumption}
\newtheorem{example}{Example}
\newtheorem{definition}{Definition}
\newtheorem{cor}{Corollary}
\newcommand{\ie}{i.e. }
\newcommand{\eg}{e.g. }
\newcommand{\cf}{cf.}
\newcommand{\R}[1]{\mathcal{R}^{#1}}
\newcommand{\PD}[1]{\mathcal{S}_{++}^{#1}}
\newcommand{\PSD}[1]{\mathcal{S}_{+}^{#1}}
\newcommand{\Null}[1]{\mathcal{N}(#1)}
\newcommand{\Range}[1]{\mbox{Im}(#1)}

\maketitle
\begin{abstract}                          
The alternating direction method of multipliers (ADMM) has emerged as a powerful technique for large-scale structured optimization. Despite many recent results on the convergence properties of ADMM, a quantitative characterization of the impact of the algorithm parameters on the convergence times of the method is still lacking. In this paper we find the optimal algorithm parameters that minimize the convergence factor of the ADMM iterates in the context of $\ell_2$-regularized minimization and constrained quadratic programming. Numerical examples show that our parameter selection rules significantly outperform existing alternatives in the literature.
\end{abstract}

\section{Introduction}
\label{sec:introduction}

The alternating direction method of multipliers  is a powerful algorithm for solving structured convex optimization problems. While the ADMM method was introduced for optimization in the 1970's, its origins can be traced back to
techniques for solving elliptic and parabolic partial difference equations developed in the 1950's~(see \cite{Boyd11} and references therein).
ADMM enjoys the
strong convergence properties of the method of multipliers and the decomposability property of dual ascent, and is particularly useful for solving optimization problems that are too large to be handled by generic optimization solvers. The method has found a large number of
applications in diverse areas such as compressed sensing~\cite{Yang:2011}, regularized estimation~\cite{bo2012},  image processing~\cite{image010}, machine learning~\cite{Forero010}, and resource allocation in wireless networks~\cite{Joshi:12}. This broad range of applications has triggered a strong recent interest in developing a better understanding of the theoretical properties of ADMM\cite{deng12,luo12,boley:13}.

Mathematical decomposition is a classical approach for parallelizing numerical optimization algorithms. If the decision problem has a favorable structure, decomposition techniques such as primal and dual decomposition allow to distribute the computations on multiple processors\cite{Las:70,BeT:89}. The processors are coordinated towards optimality by solving a suitable master problem, typically using gradient or subgradient techniques. If problem parameters such as Lipschitz constants and convexity parameters of the cost function are available, the optimal step-size parameters and associated convergence rates are well-known (\emph{e.g.}, \cite{Nesterov03}). A drawback of the gradient method is that it is sensitive to the choice of the step-size, even to the point where poor parameter selection can lead to algorithm divergence.
In contrast, the ADMM technique is surprisingly robust to poorly selected algorithm parameters: under mild conditions, the method is guaranteed to converge for all positive values of its single parameter.
Recently, an intense research effort has been devoted to establishing the rate of convergence of the ADMM method. It is now known that if the objective functions are strongly convex and have Lipschitz-continuous gradients, then the iterates produced by the ADMM algorithm converge linearly to the optimum in a certain distance metric~\eg\cite{deng12}. The application of ADMM to quadratic problems was considered in~\cite{boley:13}  and it was conjectured that the iterates converge linearly in the neighborhood of the optimal solution. It is important to stress that even when the ADMM method has linear convergence \emph{rate}, the number of iterations ensuring a desired accuracy, \ie the convergence \emph{time}, is heavily affected by the choice of the algorithm parameter. We will show that a poor parameter selection
can result in arbitrarily large convergence times for the ADMM algorithm.

The aim of the present paper is to contribute to the understanding of the convergence properties of the ADMM method. Specifically, we derive the algorithm parameters that minimize the convergence factor of the ADMM iterations for two classes of quadratic optimization problems: $\ell_2$-regularized quadratic minimization and quadratic programming with linear inequality constraints. In both cases,  we establish linear convergence rates and develop techniques to minimize the convergence factors of the ADMM iterates. These techniques allow us to give explicit expressions for the optimal algorithm parameters and the associated convergence factors. We also study over-relaxed ADMM iterations and demonstrate how to jointly choose the ADMM parameter and the over-relaxation parameter to improve the convergence times even further.
We have chosen to focus on quadratic problems, since they allow for analytical tractability, yet have vast applications in estimation~\cite{Falcao1995}, multi-agent systems~\cite{nedic10} and control\cite{DSB:2013}. Furthermore,  many complex problems can be reformulated as or approximated by QPs~\cite{SBV:2004}, and optimal ADMM parameters for QP's can be used as a benchmark for more complex ADMM sub-problems~\eg $\ell_1$-regularized problems~\cite{Boyd11}.
To the best of our knowledge, this is one of the first works that addresses the problem of optimal parameter selection for ADMM. A few recent papers have focused on the optimal parameter selection of ADMM algorithm for some variations of distributed convex programming subject to linear equality constraints \eg\cite{TGS:13,SLY:13}.

The paper is organized as follows. In Section~\ref{sec:backgroud}, we derive some preliminary results on fixed-point iterations and review the necessary background on the ADMM method. Section~\ref{sec:l2} studies $\ell_2$-regularized quadratic programming and gives explicit expressions for the jointly optimal step-size and acceleration parameter that minimize the convergence factor. We then shift our focus to the quadratic programming with linear inequality constraints  and derive the optimal step-sizes for such problems in Section~\ref{sec:qp}. We  also consider two acceleration techniques and discuss inexpensive ways to improve the speed of convergence. Our results are illustrated through numerical examples in Section~\ref{sec:qp_evaluation}. In Section~\ref{sec:qp_evaluation} we perform an extensive Model Predictive Control (MPC) case study and evaluate the performance of ADMM with the proposed parameter selection rules. A comparison with an accelerated ADMM method from the literature is also performed. Final remarks and future directions conclude the paper.

\subsection{Notation}
We denote the set of real numbers with $\R{}$ and define the set of positive (nonnegative) real numbers as $\R{}_{++}$ ($\R{}_{+}$). Let ${\mathcal S^n}$ be the set of real symmetric matrices of dimension $n\times n$. The set of positive definite (semi-definite) $n\times n$ matrices is denoted by $\PD{n}$ ($\PSD{n}$). With $I$ and $I_m$,  we symbolize the identity matrix and the identity matrix of a dimension $m\times m$, respectively.

Given a matrix $A\in \R{n\times m}$, let $\Null{A}\triangleq\{x\in\R{m} \vert \; Ax=0\}$ be the null-space of $A$ and denote the range space of $A$ by $\Range{A}\triangleq\{y\in\R{n} \vert \; y=Ax,\;x\in\R{m}\}$. We say the nullity of $A$ is $0$ (of zero dimensional) when $\Null{A}$ only contains $0$. The transpose of $A$ is represented by $A^\top$ and for $A$ with full-column rank we define $A^\dagger \triangleq (A^\top A)^{-1}A^\top$ as the pseudo-inverse of $A$. Given a subspace $\mathcal{X}\subseteq \R{n}$, $\Pi_{\mathcal{X}}\in\R{n\times n}$ denotes the orthogonal projector onto $\mathcal{X}$, while $\mathcal{X}^\bot$ denotes the orthogonal complement of $\mathcal{X}$.

For a square matrix $A$ with an eigenvalue $\lambda$ we call the space spanned by all the eigenvectors corresponding to the eigenvalue $\lambda$ the $\lambda$-eigenspace of $A$. The $i$-th smallest in modulus eigenvalue is indicated by $\lambda_i(\cdot)$. The spectral radius of a matrix $A$ is denoted by $r(A)$. The vector (matrix) $p$-norm is denoted by $\Vert \cdot \Vert_p$ and $\Vert \cdot \Vert = \Vert \cdot \Vert_2$ is the Euclidean (spectral) norm of its vector (matrix) argument. Given a subspace $\mathcal{X}\subseteq\R{n}$ and a matrix $A\in\R{n\times n}$, denote $\|A\|_{\mathcal{X}} = \max_{x\in\mathcal{X}}\dfrac{\Vert Ax \Vert}{\Vert x\Vert}$ as the spectral norm of $A$ restricted to the subspace $\mathcal{X}$.

Given $z\in \R{n}$, the diagonal matrix $Z\in \R{n\times n}$ with $Z_{ii} = z_i$
and $Z_{ij}=0$ for $j\neq i$
is denoted by $Z=\mbox{diag}(z)$. Moreover, $z\geq0$ denotes the element-wise inequality, $\vert z \vert$ corresponds to the element-wise absolute value of $z$, and $\mathcal{I}_{+}(z)$ is the indicator function of the positive orthant defined as $\mathcal{I}_{+}(z) = 0$ for $z\geq0$ and $\mathcal{I}_{+}(z) = +\infty$ otherwise.

 Consider a sequence $\{x^k\}$ converging to a fixed-point $x^\star\in\R{n}$. The \emph{convergence factor} of the converging sequence is defined as
\begin{align}\label{eqn:convergence_factor_def}
\zeta&\triangleq\, \underset{k:\,x^k \neq x^\star}{\mbox{sup}} \dfrac{\Vert x^{k+1}-x^\star\Vert}{\Vert x^{k}-x^\star\Vert}.
\end{align}
The sequence $\{x^k\}$ is said to converge Q-sublinearly if $\zeta = 1$, Q-linearly if $\zeta^k\in (0,1)$, and Q-superlinearly if $\zeta = 0$. Moreover, we say that convergence is R-linear if there is a nonnegative scalar sequence $\{\nu_k\}$ such that $\Vert x^k-x^\star \Vert\leq \nu_k$ for all $k$ and $\{\nu_k\}$ converges Q-linearly to $0$~\cite{JNW:2006}~\footnote{The letters Q and R stand for quotient and root, respectively.}. In this paper, we omit the letter Q while referring the convergence rate.

Given an initial condition $x^0$ such that $\Vert x^{0}-x^\star\Vert\leq \sigma $, we define the \emph{$\varepsilon$-solution time} $\pi_{\varepsilon}$ as the smallest iteration count to ensure that $\Vert x_k\Vert \leq \varepsilon$ holds for all $k\geq \pi_{\varepsilon}$. For linearly converging sequences with $\zeta\in(0,1)$ the $\varepsilon$-solution time is given by $\pi_\varepsilon  \triangleq\, \dfrac{\log(\sigma)-\log({\varepsilon})} {-\log(\zeta)}$.
If the $0$-solution time is finite for all $x^0$, we say that the sequence converges in finite time. As for linearly converging sequences $\zeta < 1$, the $\varepsilon$-solution time $\pi_\varepsilon$ is reduced by minimizing $\zeta$.
%
\section{Background and preliminaries}
\label{sec:backgroud}
This section presents preliminary results on fixed-point iterations and the ADMM method.

\subsection{Fixed-point iterations}

Consider the following iterative process
\begin{align}
\label{eqn:fixed_point_iterates}
x^{k+1} = T x^{k},
\end{align}
where $x^k\in \R{n}$ and $T \in \mathcal{S}^{n\times n}$.
Assume $T$ has $m<n$ eigenvalues at $1$ and let $V\in\R{n\times m}$ be a matrix whose columns span the $1$-eigenspace of $T$ so that $TV=V$.

Next we determine the properties of $T$ such that, for any given starting point $x^0$, the iteration in~\eqref{eqn:fixed_point_iterates} converges to a fixed-point that is the projection of the $x^0$ into the $1$-eigenspace of $T$, i.e.
\begin{align}
\label{eqn:fixed_point_limit}
x^\star\triangleq\lim_{k\rightarrow \infty} x^k = \lim_{k\rightarrow \infty} T^k x^0 = \Pi_{\Range{V}}x^0. 
\end{align}

\begin{proposition}
\label{prop:1}
The iterations~\eqref{eqn:fixed_point_iterates} converge to a fixed-point in $\Range{V}$ if and only if
\begin{align}
\label{eqn:convergence_condition}
r\left(T-\Pi_{\Range{V}} \right)<1.
\end{align}
\end{proposition}
\begin{IEEEproof}
The result is an extension of~\cite[Theorem~1]{XiB:04} for the case of $1$-eigenspace of $T$ with dimension $m>1$. The proof is similar to this citation and is therefore omitted.
\end{IEEEproof}

Proposition~\ref{prop:1} shows that when $T\in \mathcal{S}^n$, the fixed-point iteration~\eqref{eqn:fixed_point_iterates} is guaranteed to converge to a point given by (\ref{eqn:fixed_point_limit}) if all the non-unitary eigenvalues of $T$ have magnitudes strictly smaller than 1. 
From~\eqref{eqn:fixed_point_iterates} one sees that
\begin{align*}
x^{k+1}-x^\star &= \left(T-\Pi_{\Range{V}}\right) x^k = \left(T-\Pi_{\Range{V}}\right) (x^k-x^\star)
\end{align*}
Hence, the convergence factor of~\eqref{eqn:fixed_point_iterates} is the modulus of the largest non-unit eigenvalue of $T$. 
\subsection{The ADMM method}
The ADMM algorithm solves problems of the form
\begin{align}
	\begin{array}[c]{ll}
	\mbox{minimize} & f(x)+g(z)\\
	\mbox{subject to} & Ax+Bz=c
	\end{array} \label{eqn:admm_standard_form}
\end{align}
where $f$ and $g$ are convex functions, $x\in {\mathcal R}^n$, $z\in {\mathcal R}^m$, $A\in {\mathcal R}^{p\times n}$, $B\in {\mathcal R}^{p\times m}$ and $c\in {\mathcal R}^p$; see~\cite{Boyd11} for a detailed review.

Relevant examples that appear in this form are, \eg regularized estimation, where $f$ is the estimator loss and $g$ is the regularization term,  and various networked optimization problems, \emph{e.g.}~\cite{Italian,Boyd11}.
The method is based on the \emph{augmented Lagrangian}
\begin{align*}
	L_{\rho}(x,z,\mu) &= f(x)+g(z) + \dfrac{\rho}{2}\Vert Ax+Bz-c\Vert_2^2 + \mu^{T}(Ax+Bz-c),
\end{align*}
and performs sequential minimization of the $x$ and $z$ variables followed by a dual variable update:
\begin{align}
	x^{k+1} &= \underset{x}{\operatorname{argmin}}\, L_{\rho}(x,z^{k}, \mu^{k}), \nonumber\\
	z^{k+1} &= \underset{z}{\operatorname{argmin}}\, L_{\rho}(x^{k+1}, z, \mu^k), \label{eqn:admm_iterations}\\
	\mu^{k+1} &= \mu^{k} + \rho(Ax^{k+1}+Bz^{k+1}-c), \nonumber
\end{align}
for some arbitrary $x^0 \in  {\mathcal R}^n$, $z^0\in {\mathcal R}^m$, and $\mu^0\in {\mathcal R}^p$.
It is often convenient to express the iterations in terms of the scaled dual variable $u=\mu/\rho$:
\begin{align}
\label{eqn:admm_scaled}
\begin{array}[c]{ll}
	x^{k+1} &= \underset{x}{\operatorname{argmin}} \left\{ f(x)+ \dfrac{\rho}{2}\Vert Ax+Bz^k-c+u^k\Vert_2^2\right\}, \\
	z^{k+1} &= \underset{z}{\operatorname{argmin}} \left\{g(z) + \dfrac{\rho}{2}\Vert Ax^{k+1}+Bz-c+u^{k}\Vert_2^2\right\}, \\
	u^{k+1} &= u^{k} + Ax^{k+1}+Bz^{k+1}-c.
\end{array}
\end{align}
ADMM is particularly useful when the $x$- and $z$-minimizations can be carried out efficiently, for example when they admit closed-form expressions. Examples of such problems include linear and quadratic programming, basis pursuit, $\ell_1$-regularized minimization, and model fitting problems to name a few (see~\cite{Boyd11} for a complete discussion). One advantage of the ADMM method is that there is only a single algorithm parameter, $\rho$, and under rather mild conditions, the method can be shown to converge for all values of the parameter; see \cite{Boyd11, JEK:12} and references therein. As discussed in the introduction, this contrasts the gradient method  whose iterates diverge if the step-size parameter is chosen too large. However, $\rho$ has a direct impact on the convergence factor of the algorithm, and inadequate tuning of this parameter can render the method slow.
The convergence of ADMM is often characterized in terms of the residuals
\begin{align}
r^{k+1} &= Ax^{k+1}+B z^{k+1}-c,\label{eq:primal_res}\\
s^{k+1} &= \rho A^\top B (z^{k+1}-z^k),\label{eq:dual_res}
\end{align}
termed the \emph{primal} and \emph{dual} residuals, respectively~\cite{Boyd11}. One approach for improving the convergence properties of the algorithm is to  also account for past iterates when computing the next ones. This technique is called \emph{relaxation} and  amounts to replacing $A x^{k+1}$ with $h^{k+1} = \alpha^k A x^{k+1}- (1-\alpha^k) (Bz^k -c)$ in the $z$- and $u$-updates \cite{Boyd11}, yielding
\begin{equation}
\label{eqn:admm_relaxed}
\begin{aligned}
%
	z^{k+1} &= \underset{z}{\operatorname{argmin}} \,\left\{g(z) + \dfrac{\rho}{2}\left\Vert h^{k+1} +Bz-c+u^{k}\right\Vert_2^2\right\},  \\
	u^{k+1} &= u^{k} + h^{k+1}+ Bz^{k+1}-c.
\end{aligned}
\end{equation}
The parameter $\alpha^k\in (0,2)$ is called the \emph{relaxation parameter}. Note that letting $\alpha^k=1$ for all $k$ recovers the original ADMM iterations~\eqref{eqn:admm_scaled}. Empirical studies show that over-relaxation, \emph{i.e.} letting ${\alpha^k>1}$, is often advantageous and the guideline $\alpha^k\in [1.5, 1.8]$ has been proposed~\cite{Eckstein:1994}.


In the rest of this paper, we will consider the traditional ADMM iterations~\eqref{eqn:admm_iterations} and the relaxed version~\eqref{eqn:admm_relaxed} for different classes of quadratic problems, and derive explicit expressions for the step-size $\rho$ and the relaxation parameter $\alpha$ that minimize the convergence factors.
\section{ Optimal convergence factor for $\ell_2$-regularized quadratic minimization}
\label{sec:l2}
Regularized estimation problems
\begin{align*}
	\begin{array}[c]{ll} \mbox{minimize} & f(x) + \dfrac{\delta}{2}\Vert x \Vert_p^q \end{array}
\end{align*}
where $\delta>0$ are abound in statistics, machine learning, and control. In particular, $\ell_1$-regularized estimation where $f(x)$ is quadratic and $p=q=1$, and \emph{sum of norms} regularization where $f(x)$ is quadratic, $p=2$, and $q=1$, have recently received significant attention  \cite{ohlsson2010segmentation}.
In this section we will focus on $\ell_2$-regularized estimation, where $f(x)$ is quadratic and $p=q=2$, \ie
\begin{align}
\label{eqn:L2_formulation}
\begin{array}{ll}
\mbox{minimize} & \dfrac{1}{2} x^\top Q x + q^\top x + \dfrac{\delta}{2} \Vert z\Vert_2^2 \\
\mbox{subject to}& x-z=0,
\end{array}
\end{align}
for $Q \in \PD{n}$, $x,q,z \in \R{n}$ and constant regularization parameter $\delta \in \mathcal{R}_+$. While these problems can be solved explicitly and do not motivate the ADMM machinery per se, they provide insight into the step-size selection for ADMM and allow us to compare the performance of an optimally tuned ADMM to direct alternatives
(see Section~\ref{sec:qp_evaluation}).

\subsection{Standard ADMM iterations}
The standard ADMM iterations are given by
\begin{align}
\label{eqn:ADMM_L2_iterations}
 \begin{array}{l}
 x^{k+1}=(Q+\rho I)^{-1}(\rho z^k-\mu^k-q),\\
 z^{k+1} = \dfrac{\mu^k+\rho x^{k+1}}{\delta+\rho},\\
 \mu^{k+1}= \mu^{k}+\rho (x^{k+1}-z^{k+1}).
 \end{array}
 \end{align}
The $z$-update implies that  $\mu^k=(\delta+\rho)z^{k+1}-\rho x^{k+1}$, so the $\mu$-update can be re-written as
\begin{align}
\nonumber
\mu^{k+1}=(\delta + \rho)z^{k+1}-\rho x^{k+1} + \rho(x^{k+1}-z^{k+1})= \delta z^{k+1}.
\end{align}
Hence, to study the convergence of (\ref{eqn:ADMM_L2_iterations}) one can investigate how the errors associated with $x^k$ or $z^k$ vanish. Inserting the $x$-update into the $z$-update and using the fact that $\mu^k=\delta z^k$,  we find
\begin{equation}
\begin{aligned}\label{eqn:ADMM_L2_matrix}
z^{k+1}&= \underbrace{\dfrac{1}{\delta+\rho}\left(\delta I + \rho (\rho-\delta)\left(Q+\rho I\right)^{-1}\right)}_E z^k -\dfrac{\rho}{\delta+\rho}(Q+\rho I)^{-1}q.
\end{aligned}
\end{equation}
Let $z^{\star}$ be a fixed-point 
of \eqref{eqn:ADMM_L2_matrix}, i.e. $z^{\star}= E z^\star -\dfrac{\rho(Q+\rho I)^{-1}}{\delta+\rho}q$. The dual error $e^{k+1}\triangleq z^{k+1}-z^\star$ then evolves as
\begin{align}
	e^{k+1} &= E e^{k}.\label{eqn:ADMM_L2_error}
\end{align}
A direct analysis of the error dynamics (\ref{eqn:ADMM_L2_error}) allows us to characterize the convergence of (\ref{eqn:ADMM_L2_iterations}):
\begin{theorem}
\label{thm:L2:standard}
For all values of the step-size $\rho>0$ and regularization parameter $\delta >0$, both $x^k$ and $z^k$ in the ADMM iterations~(\ref{eqn:ADMM_L2_iterations}) converge to $x^\star=z^\star$, the solution of optimization problem~\eqref{eqn:L2_formulation}. Moreover, $z^{k+1}-z^\star$ converges at linear rate $\zeta \in (0,1)$ for all $k\geq 0$. The pair of the optimal constant step-size $\rho^{\star}$ and convergence factor $\zeta^\star$ are given as
\begin{equation}
\begin{aligned}
\label{eqn:ADMM_L2_optimal_step-size}
\rho^{\star} &=
\begin{cases}
\sqrt{\delta \lambda_1(Q)}&\quad \mbox{if} \; \delta< \lambda_1(Q),\\
\sqrt{\delta \lambda_n(Q)}&\quad \mbox{if} \; \delta> \lambda_n(Q),\\
\delta  &\quad\mbox{otherwise}.
\end{cases}\quad
 \zeta^{\star} &=
\begin{cases}
\left(1+\dfrac{\delta+\lambda_1(Q)}{2\sqrt{\delta \lambda_1(Q)}}\right)^{-1}& \mbox{if} \; \delta< \lambda_1(Q),\\
\left(1+\dfrac{\delta+\lambda_n(Q)}{2\sqrt{\delta \lambda_n(Q)}}\right)^{-1}& \mbox{if} \; \delta> \lambda_n(Q),\\
\dfrac{1}{2} & \mbox{otherwise}.
\end{cases}
\end{aligned}
\end{equation}
\begin{IEEEproof}
See appendix for this and the rest of the proofs.
\end{IEEEproof}
\end{theorem}
%
\begin{cor}
\label{cor:L2:standard} Consider  the error dynamics described by (\ref{eqn:ADMM_L2_error}) and $E$ in \eqref{eqn:ADMM_L2_matrix}.
For $\rho=\delta$,
\begin{align*}
	\lambda_i(E)=1/2,\qquad i=1, \dots, n,
\end{align*}
and the convergence factor of the error dynamics  \eqref{eqn:ADMM_L2_error} is independent of $Q$.
\end{cor}
\begin{rem}
Note that the convergence factors in Theorem~\ref{thm:L2:standard} and Corollary~\ref{cor:L2:standard} are guaranteed for all initial values, and that iterates generated from specific initial values might converge even faster. Furthermore, the results focus on the dual error. For example, in Algorithm~\eqref{eqn:ADMM_L2_iterations} with $\rho = \delta$ and initial condition $z^0=0$, $\mu^0 = 0$, the $x$-iterates converge in one iteration since $x^1 = -(Q+\delta I)^{-1}q=x^\star$. However, the constraint in~\eqref{eqn:L2_formulation} is not satisfied and a straightforward calculation shows that  $e^{k+1}=1/2 e^k$. Thus, although $x^k=x^\star$ for $k\geq 1$, the dual residual $\Vert e^k\Vert=\Vert z^k-z^\star\Vert$ decays linearly with a factor of $1/2$.
\end{rem}
\begin{rem}
The analysis above also applies to the more general case with cost function $\dfrac{1}{2} \bar{x}^\top \bar{Q} \bar{x} + \bar{q}^\top \bar{x} + \dfrac{\delta}{2} \bar{z}^\top \bar{P} \bar{z}$ where $\bar{P} \in \PD{n}$. A change of variables $z=\bar{P}^{1/2}\bar{z}$ is then applied to transform the problem into the form (\ref{eqn:L2_formulation}) with $x=\bar{P}^{1/2} \bar{x}$, $q= \bar{P}^{-1/2}\bar{q}$, and $Q= \bar{P}^{-1/2}\bar{Q}\bar{P}^{-1/2}$.
\end{rem}
\subsection{Over-relaxed ADMM iterations}
The over-relaxed ADMM iterations for (\ref{eqn:L2_formulation}) can be found by replacing $x^{k+1}$ by $\alpha x^{k+1} + (1-\alpha)z^k$ in the $z-$ and $\mu$-updates of \eqref{eqn:ADMM_L2_iterations}.
The resulting iterations take the form
\begin{align}
\label{eqn:ADMM_L2_iterations_relaxation}
\begin{array}[c]{ll}
 x^{k+1}&=(Q+\rho I)^{-1}(\rho z^k-\mu^k-q),\\
 z^{k+1} &= \dfrac{\mu^k+\rho (\alpha x^{k+1}+ (1-\alpha) z^k)}{\delta+\rho},\\
 \mu^{k+1} &= \mu^{k}+\rho \left(\alpha (x^{k+1} - z^{k+1}) + (1-\alpha)\left(z^k - z^{k+1}\right)\right).
\end{array}
\end{align}
The next result demonstrates that in a certain range of $\alpha$ it is possible to obtain a guaranteed improvement of the convergence factor compared to the classical iterations~\eqref{eqn:ADMM_L2_iterations}.
\begin{theorem}\label{thm:L2:relaxed}
Consider the $\ell_2$-regularized quadratic minimization problem~\eqref{eqn:L2_formulation} and its associated over-relaxed ADMM iterations~\eqref{eqn:ADMM_L2_iterations_relaxation}.
For all positive step-sizes $\rho>0$ and all relaxation parameters $\alpha\in (0, 2\underset{i}{\min}\{ (\lambda_i(Q)+\rho)(\rho+\delta)/(\rho\delta + \rho \lambda_i(Q))\})$, the iterates $x^k$ and $z^k$ converge to the solution of~\eqref{eqn:L2_formulation}.
Moreover, the dual variable converges at linear rate $\Vert z^{k+1} - z^\star\Vert \leq \zeta_R \Vert z^k - z^\star\Vert$ and the convergence factor $\zeta_R <1$ is strictly smaller than that of the classical ADMM algorithm~\eqref{eqn:ADMM_L2_iterations} if $1<\alpha<2\underset{i}{\min}\{ (\lambda_i(Q)+\rho)(\rho+\delta)/(\rho\delta + \rho \lambda_i(Q))\}$
The jointly optimal step-size, relaxation parameter, and the convergence factor $(\rho^\star, \alpha^\star,\zeta_R^\star)$ are given by
 \begin{align}
\label{eqn:ADMM_L2_Relaxation_optimal}
\rho^\star = \delta, \quad \alpha^\star = 2,\quad \zeta^\star_R = 0.
\end{align}
With these parameters, the ADMM iterations converge in one iteration.

\end{theorem}

\begin{rem}
The upper bound on $\alpha$ which ensures faster convergence of the over-relaxed ADMM iterations~\eqref{eqn:ADMM_L2_iterations_relaxation} compared to~\eqref{eqn:ADMM_L2_iterations} depends on the eigenvalues of $Q$, $\lambda_i(Q)$, which might be unknown. However, since $(\rho+\delta)(\rho+\lambda_i(Q))> \rho(\lambda_i(Q)+\delta)$ the over-relaxed iterations are guaranteed to converge faster for all $\alpha \in (1,2]$, independently of $Q$.
\end{rem}

\section{Optimal convergence factor for quadratic programming}
\label{sec:qp}
In this section, we consider a quadratic programming (QP) problem of the form
\begin{align}
\label{eqn:Quadratic_problem}
	\begin{array}[c]{ll}
	\mbox{minimize} & \dfrac{1}{2} x^\top Q x+ q^\top x\\
	\mbox{subject to}& Ax \leq c
	\end{array}
\end{align}
where $Q\in \PD{n}$, $q \in \R{n}$, $A\in \mathcal{R}^{m\times n}$ is full rank and $c\in \mathcal{R}^{m}$.

\subsection{Standard ADMM iterations}
The QP-problem (\ref{eqn:Quadratic_problem}) can be put on ADMM standard form  (\ref{eqn:admm_standard_form}) by introducing a slack vector $z$ and putting an infinite penalty on negative components of $z$, \emph{i.e.}
\begin{align}
\label{eqn:Quadratic_problem_1}
	\begin{array}[c]{ll}
	\mbox{minimize} & \dfrac{1}{2} x^\top Q x+ q^\top x+ \mathcal{I}_{+}(z)\\
	\mbox{subject to}& Ax - c + z = 0.
	\end{array}
\end{align}
The associated \emph{augmented Lagrangian} is
\begin{align*}
L_\rho(x,z,u) = \dfrac{1}{2} x^\top Q x+ q^\top x+ \mathcal{I}_{+}(z) + \dfrac{\rho}{2} \Vert Ax - c + z + u \Vert^2_2,
\end{align*}
where $u=\mu/\rho$, which leads to the scaled ADMM iterations
\begin{align}
\begin{array}[c]{ll}
	x^{k+1} &= -(Q+\rho A^\top A)^{-1} [q+\rho A^\top(z^k + u^k - c)], \\
	z^{k+1} &= \mbox{max}\{0,-A x^{k+1}-u^{k}+c\}, \\
	u^{k+1} &= u^{k} + A x^{k+1}-c+z^{k+1}.
\end{array}
\label{eqn:Quadratic_admm_iterations}
\end{align}
To study the convergence of~\eqref{eqn:Quadratic_admm_iterations} we rewrite it in an equivalent form with linear time-varying matrix operators. To this end, we introduce a vector of indicator variables $d^k\in \{0,1\}^{n}$  such that $d^k_i = 0$ if $u^k_i=0$ and $d^k_i=1$ if $u^k_i\neq 0$.  From the $z$- and $u$- updates in~\eqref{eqn:Quadratic_admm_iterations}, one observes that $z_i^k\neq 0\rightarrow u_i^k=0$, \emph{i.e.} $u_i^k\neq 0\rightarrow z_i^k=0$.
Hence, $d^k_i =1$ means that at the current iterate, the slack variable $z_i$ in~\eqref{eqn:Quadratic_problem_1} equals zero; i.e., the $i^{\rm th}$ inequality constraint in~\eqref{eqn:Quadratic_problem} is active. We also introduce the variable vector $v^k\triangleq z^k+u^k$ and let $D^k=\mbox{diag}(d^k)$ so that $D^k v^k= u^k$ and $(I-D^k)v^k=z^k$. Now, the second and third steps of~\eqref{eqn:Quadratic_admm_iterations} imply that $v^{k+1} = \left\vert Ax^{k+1}+u^k -c\right\vert = F^{k+1}(A x^{k+1}+ D^k v^k-c)$ where $F^{k+1}\triangleq \mbox{diag}\left ( \text{sign}(A x^{k+1}+ D^k v^k-c) \right )$ and $\text{sign}(\cdot)$ returns the signs of the elements of its vector argument.
 Hence,~\eqref{eqn:Quadratic_admm_iterations} becomes
\begin{align}
\begin{array}[c]{ll}
	x^{k+1} &= -(Q+\rho A^\top A)^{-1} [q+\rho A^\top (v^k-c)], \\
	v^{k+1} &= \left\vert A x^{k+1}+ D^k v^k - c\right\vert = F^{k+1}(A x^{k+1}+ D^k v^k-c), \\
	D^{k+1} &= \dfrac{1}{2}(I+F^{k+1}),
\end{array}
\label{eqn:Quadratic_admm_iterations_reformed}
\end{align}
where the $D^{k+1}$-update follows from the observation that
\begin{align}
\nonumber
(D_{ii}^{k+1},\, F_{ii}^{k+1}) =\left\{
\begin{array}[l]{lll}
\hspace{-5pt}(0, \,-1) & \hspace{-3pt}\mbox{if} & \hspace{-2pt}v_i^{k+1} = -(A x_i^{k+1}+u_i^k - c) \\
\hspace{-5pt}(1,\, 1) & \hspace{-3pt}\mbox{if} & \hspace{-2pt}v_i^{k+1} = A x_i^{k+1}+ u_i^k- c
\end{array}\right.
\end{align}
%
 Since the $v^k$-iterations will be central in our analysis, we will develop them further. Inserting the expression for $x^{k+1}$ from the first equation of~\eqref{eqn:Quadratic_admm_iterations_reformed} into the second, we find
\begin{equation}
\begin{aligned}
\label{eqn:v_recurrence}
v^{k+1} &= F^{k+1}\Big(\left(  D^k - A (Q/\rho+ A^\top A)^{-1} A^\top \right) v^k\Big)  - F^{k+1}\Big(A (Q+\rho A^\top A)^{-1}(q-\rho A^\top c) + c\Big).
\end{aligned}
\end{equation}
Noting that $D^k=\dfrac{1}{2}(I+F^k)$ and introducing
\begin{align}
M &\triangleq  A (Q/\rho+ A^\top A)^{-1} A^\top,
\end{align}
we obtain
\begin{equation}
\begin{aligned}
\label{eqn:QP_Fv_sequence}
F^{k+1}v^{k+1} - F^k v^k &= \left(  \dfrac{I}{2} - M \right) (v^k - v^{k-1}) + \dfrac{1}{2}\left( F^k v^k - F^{k-1} v^{k-1}\right).
\end{aligned}
\end{equation}
 We now relate $v^k$ and $F^k v^k$ to the primal and dual residuals, $r^k$ and $s^k$, defined in (\ref{eq:primal_res}) and (\ref{eq:dual_res}):
\begin{proposition}\label{prop:w_residuals}
Consider $r^{k}$ and $s^{k}$ the primal and dual residuals of the QP-ADMM algorithm~\eqref{eqn:Quadratic_admm_iterations} and auxiliary variables $v^k$ and $F^k$. The following relations hold
\begin{align}
&F^{k+1}v^{k+1} - F^kv^{k} = r^{k+1} - \dfrac{1}{\rho}Rs^{k+1} - \Pi_{\Null{A^\top}}(z^{k+1} - z^k), \label{eqn:w_minus_residual}\\
&v^{k+1} - v^{k} = r^{k+1} + \dfrac{1}{\rho}Rs^{k+1} + \Pi_{\Null{A^\top}}(z^{k+1} - z^k),\label{eqn:w_plus_residual}
\end{align}
\begin{align}\label{eqn:r_to_fv}
    & \Vert r^{k+1} \Vert \leq \Vert F^{k+1} v^{k+1} - F^{k} v^{k} \Vert,\\
    & \Vert s^{k+1} \Vert \leq \rho \Vert A\Vert \Vert F^{k+1} v^{k+1} - F^{k} v^{k} \Vert. \label{eqn:s_to_fv}
\end{align}
where
\begin{enumerate}
 \item[(i)] $R= A(A^\top A)^{-1}$ and $\Pi_{\Null{A^\top}} = I-A(A^\top A)^{-1}A^\top$, if $A$ has full column-rank;
\item[(ii)] $R=(AA^\top)^{-1}A$ and $\Pi_{\Null{A^\top}} = 0$, if $A$ has full row-rank;
\item[(iii)] $R = A^{-1}$ and $\Pi_{\Null{A^\top}} = 0$, if $A$ is invertible.
\end{enumerate}
\end{proposition}
The next theorem guarantees that~\eqref{eqn:QP_Fv_sequence} convergence linearly to zero in the auxiliary residuals~\eqref{eqn:w_minus_residual} which implies R-linear convergence of the ADMM algorithm~\eqref{eqn:Quadratic_admm_iterations} in terms of the primal and dual residuals. The optimal step-size $\rho^\star$ and the smallest achievable convergence factor are characterized immediately afterwards.


\begin{theorem}\label{thm:QP:linear_rate}
Consider the QP~\eqref{eqn:Quadratic_problem} and the corresponding ADMM iterations~\eqref{eqn:Quadratic_admm_iterations}. For all values of the step-size $\rho\in\R{}_{++}$ the residual $ F^{k+1}v^{k+1} - F^{k}v^{k}$ converges to zero at linear rate. Furthermore, $ r^k $ and $ s^k $, the primal and dual residuals of~\eqref{eqn:Quadratic_admm_iterations}, converge R-linearly  to zero.
\end{theorem}
\begin{theorem}\label{thm:QP_optimal_factor}
Consider the {QP~\eqref{eqn:Quadratic_problem} and the corresponding ADMM iterations~\eqref{eqn:Quadratic_admm_iterations}}. If the constraint matrix $A$ is either full row-rank or invertible then the optimal step-size and convergence factor for the $F^{k+1}v^{k+1}-F^k v^k$ residuals are
\begin{equation}
\begin{aligned}\label{eqn:QP_optimal_factor}
\rho^\star &= \left(\sqrt{\lambda_1(A Q^{-1} A^\top) \lambda_n(A Q^{-1} A^\top)}\right)^{-1},\\
\zeta^\star &= \dfrac{\lambda_n(A Q^{-1} A^\top)}{\lambda_n(A Q^{-1} A^\top) + \sqrt{\lambda_1(A Q^{-1} A^\top)\lambda_n(A Q^{-1} A^\top)}}.
\end{aligned}
\end{equation}
\end{theorem}

%

Although the convergence result of Theorem~\ref{thm:QP:linear_rate} holds for all QPs of the form~\eqref{eqn:Quadratic_problem}, optimality of the step-size choice proposed in Theorem~\ref{thm:QP_optimal_factor} is only established for problems where the constraint matrix $A$ has full row-rank or it is invertible.
However, as shown next, the convergence factor can be arbitrarily close to $1$ when rows of $A$ are linearly dependent.
\begin{theorem}\label{thm:QP_slow}
Define variables
\begin{align*}
\epsilon_k &\triangleq \dfrac{\|M(v^k - v^{k-1})\|}{\|F^kv^k - F^{k-1} v^{k-1}\|},\quad \quad \delta_k\triangleq\dfrac{ \|D^kv^k - D^{k-1} v^{k-1}\| }{\|F^kv^k - F^{k-1} v^{k-1}\|},\\
\tilde{\zeta}(\rho) &\triangleq  \max_{i:\;\lambda_i(AQ^{-1}A^\top) > 0}\left\{ \left\vert \dfrac{\rho \lambda_i(AQ^{-1}A^\top)}{1+\rho\lambda_i(A Q^{-1} A^\top)} - \dfrac{1}{2}\right\vert + \dfrac{1}{2}\right\},
\end{align*}
and $\underline\zeta^k\triangleq\vert \delta_k - \epsilon_k \vert$.

The convergence factor $\zeta$ of the residual $F^{k+1}v^{k+1} - F^kv^k$ is lower bounded by
\begin{equation}\label{eq:QP_lower_bound}
\underline{\zeta}\triangleq \max_k\; \underline\zeta^k < 1.
\end{equation}
Furthermore, given an arbitrarily small $\xi\in(0,\, \frac{1}{2})$ and $\rho>0$, we have the following results:
\begin{enumerate}
\item[(i)] the inequality $\underline{\zeta}<\tilde{\zeta}(\rho)<1$ holds for all $\delta_k\in[0,\;1]$ if and only if the nullity of $A$ is zero;
\item[(ii)] when the nullity of $A$ is nonzero and $\epsilon_k \geq 1-\xi$, it holds that $\underline{\zeta} \leq \tilde{\zeta}(\rho) + \sqrt{\dfrac{\xi}{2}}$;
\item[(iii)] when the nullity of $A$ is nonzero, $\delta_k \geq 1-\xi$, and $\|\Pi_{\Null{A^\top}}(v^k - v^{k-1}) \|/\|v^k - v^{k-1}\| \geq \sqrt{1- \xi^2/\|M\|^2}$, it follows that $\underline{\zeta}\geq 1-2\xi$.
\end{enumerate}
\end{theorem}
The previous result establishes that slow convergence can occur locally for any value of $\rho$ when the nullity of $A$ is nonzero and 
$\xi$ is small. However, as section (ii) of Theorem~\ref{thm:QP_slow} suggests, in these cases,~\eqref{eqn:QP_optimal_factor} can still work as a heuristic to reduce the convergence time if $\lambda_1(AQ^{-1}A^\top)$ is taken as the smallest nonzero eigenvalue of $AQ^{-1}A^\top$. In Section~\ref{sec:qp_evaluation}, we show numerically that this heuristic performs well with different problem setups.

\subsection{Over-relaxed ADMM iterations}
Consider the relaxation of~\eqref{eqn:Quadratic_admm_iterations} obtained by replacing
$Ax^{k+1}$ in the $z$- and $u$-updates with $\alpha A x^{k+1}- (1-\alpha) (z^k - c)$. The corresponding relaxed iterations read
\begin{align}
\begin{array}[c]{ll}
	x^{k+1} &= -(Q+\rho A^\top A)^{-1} [q+\rho A^\top(z^k + u^k - c)], \\
	z^{k+1} &= \mbox{max}\{0,-\alpha (A x^{k+1}-c)+(1-\alpha)z^k-u^{k}\}, \\
	u^{k+1} &= u^{k} + \alpha (A x^{k+1}+z^{k+1}-c)+(1-\alpha)(z^{k+1}-z^{k}).
\end{array}
\label{eqn:Quadratic_admm_iterations_relaxation}
\end{align}
In next, we study convergence and optimality properties of these iterations. We observe:
\begin{lem}\label{lem:optimal_fixed_point}
Any fixed-point of~\eqref{eqn:Quadratic_admm_iterations_relaxation} corresponds to a global optimum of~\eqref{eqn:Quadratic_problem_1}.
\end{lem}
Like the analysis of~\eqref{eqn:Quadratic_admm_iterations}, introduce $v^k = z^k+u^k$ and $d^k\in \R{n}$ with $d^k_i=0$ if $u_i^k=0$ and $d_i^k=1$ otherwise. Adding the second and the third step of~\eqref{eqn:Quadratic_admm_iterations_relaxation} yields $v^{k+1} = \left\vert \alpha (A x^{k+1}-c)-(1-\alpha)z^k + u^k \right\vert$. Moreover, $D^k = \mbox{diag}(d^k)$ satisfies $D^k v^k = u^k$ and $(I-D^k)v^k = z^k$, so~\eqref{eqn:Quadratic_admm_iterations_relaxation} can be rewritten as
\begin{align}
\begin{array}[c]{ll}
	x^{k+1} &= -(Q+\rho A^\top A)^{-1} [q+\rho A^\top (v^k-c)], \\
	v^{k+1} &= F^{k+1}\Big( \alpha \left(A x^{k+1}+D^k v^k - c \right)\Big) - F^{k+1}\Big((1-\alpha) (I-2D^k)v^k \Big), \\
	D^{k+1} &= \dfrac{1}{2}(I+F^{k+1}),
\end{array}
\label{eqn:Quadratic_admm_iterations_reformed_relaxation}
\end{align}
where  $F^{k+1}\triangleq \mbox{diag}\left ( \text{sign}\left(\alpha(A x^{k+1}+D^k v^k - c) -(1-\alpha)(I-2D^k)v^k \right) \right )$.
Defining $M \triangleq  A (Q/\rho+ A^\top A)^{-1} A^\top$ and substituting the expression for $x^{k+1}$ in~\eqref{eqn:Quadratic_admm_iterations_reformed_relaxation} into the expression for $v^{k+1}$ yields
\begin{equation}
\begin{aligned}
\label{eqn:v_recurrence_relaxation}
v^{k+1} &= F^{k+1}\Big( \left(-\alpha M + (2-\alpha)D^k - (1-\alpha)I  \right) v^k\Big) - F^{k+1}\Big( \alpha A (Q+\rho A^\top A)^{-1}(q-\rho A^\top c) + \alpha c\Big).
\end{aligned}
\end{equation}
As in the previous section, we replace $D^k$ by $\dfrac{1}{2}(I+F^k)$ in~\eqref{eqn:v_recurrence_relaxation} and form $F^{k+1}v^{k+1}-F^k v^k$:
\begin{equation}
\begin{aligned}
\label{eqn:QP_Fv_sequence_relaxation}
F^{k+1}v^{k+1} -F^k v^k &= \dfrac{\alpha}{2}\left(I-2M\right) \left(v^{k}-v^{k-1}\right) + (1-\dfrac{\alpha}{2})\left(F^{k}v^{k} -F^{k-1} v^{k-1}\right).
\end{aligned}
\end{equation}

The next theorem characterizes the convergence rate of the relaxed ADMM iterations.
\begin{theorem}
\label{thm:QP_relaxation_convergence}
Consider the QP \eqref{eqn:Quadratic_problem} and the corresponding relaxed ADMM iterations~\eqref{eqn:Quadratic_admm_iterations_relaxation}. If
\begin{align}
 \rho\in\mathcal{R}_{++}, \quad \alpha \in (0, 2],
\end{align}
then the equivalent fixed point iteration~\eqref{eqn:QP_Fv_sequence_relaxation} converges linearly in terms of $F^{k+1}v^{k+1} -F^k v^k$ residual. Moreover, $ r^k$ and $ s^k$, the primal and dual residuals of \eqref{eqn:Quadratic_admm_iterations_relaxation}, converge R-linearly to zero.
\end{theorem}
Next, we restrict our attention to the case where $A$ is either invertible or full row-rank to be able to derive the jointly optimal step-size and over-relaxation parameter, as well as an explicit expression for the associated convergence factor.
The result shows that the over-relaxed ADMM iterates can yield a significant speed up compared to the standard ADMM iterations.
\begin{theorem}
\label{thm:QP_relaxation_optimal_factor}
Consider the QP \eqref{eqn:Quadratic_problem} and the corresponding relaxed ADMM iterations~\eqref{eqn:Quadratic_admm_iterations_relaxation}. If the constraint matrix $A$ is of full row-rank or invertible then the joint optimal step-size, relaxation parameter and the convergence factor with respect to the $F^{k+1}v^{k+1} -F^k v^k $ residual are
\begin{equation}
\begin{aligned}
\label{eqn:QP_relaxation_optimal_factor}
\rho^\star &= \left(\sqrt{\lambda_1(AQ^{-1}A^\top)\; \lambda_n(AQ^{-1}A^\top)}\right)^{-1}, \quad \alpha^\star = 2,\\
   \zeta_R^\star &= \dfrac{\lambda_n(AQ^{-1}A^\top)-\sqrt{\lambda_1(AQ^{-1}A^\top)\;\lambda_n(AQ^{-1}A^\top)}}{\lambda_n(AQ^{-1}A^\top) + \sqrt{\lambda_1(AQ^{-1}A^\top)\;\lambda_n(AQ^{-1}A^\top)}}
\end{aligned}
\end{equation}
Moreover, when the iterations~\eqref{eqn:QP_Fv_sequence_relaxation} are over-relaxed;~\ie~$\alpha \in (1,2]$ their iterates have a smaller convergence factor than that of~\eqref{eqn:QP_Fv_sequence}.
\end{theorem}
\subsection{Optimal constraint preconditioning}
In this section, we consider another technique to improve the convergence of the ADMM method. The approach is based on the observation that the optimal convergence factors $\zeta^\star$ and $\zeta_R^\star$ from Theorem~\ref{thm:QP_optimal_factor} and Theorem~\ref{thm:QP_relaxation_optimal_factor}
are monotone increasing in the ratio $\lambda_n(AQ^{-1}A^\top)/\lambda_1(AQ^{-1}A^\top)$. This ratio can be decreased --without changing the complexity of the ADMM algorithm~\eqref{eqn:Quadratic_admm_iterations}-- by scaling the equality constraint in~\eqref{eqn:Quadratic_problem_1} by a diagonal matrix $L\in\PD{m}$, i,e., replacing $Ax-c+z=0$ by $L\left(Ax-c+z\right) = 0$. Let $\bar{A} \triangleq LA$, $\bar{z} \triangleq Lz$, and $\bar{c}\triangleq Lc$. The resulting scaled ADMM iterations are derived by replacing $A$, $z$, and $c$ in~\eqref{eqn:Quadratic_admm_iterations} and~\eqref{eqn:Quadratic_admm_iterations_relaxation} by the new variables $\bar{A}$, $\bar{z}$, and $\bar{c}$, respectively. Furthermore, the results of Theorem~\ref{thm:QP_optimal_factor} and Theorem~\ref{thm:QP_relaxation_optimal_factor} can be applied to the scaled ADMM iterations in terms of new variables.
Although these theorems only provide the optimal step-size parameters for the QP when the constraint matrices are invertible or have full row-rank, we use the expressions as heuristics when the constraint matrix has full column-rank. Hence, in the following we consider $\lambda_n(\bar{A} Q^{-1} \bar{A}^{\top})$ and $\lambda_1(\bar{A} Q^{-1} \bar{A}^{\top})$ to be the largest and smallest nonzero eigenvalues of $\bar{A} Q^{-1} \bar{A}^{\top} = LAQ^{-1}A^\top L$, respectively and minimize the ratio $\lambda_n/{\lambda_1}$ in order to minimize the convergence factors $\zeta^\star$ and $\zeta_R^\star$. A similar problem was also studied in~\cite{GIS:14,GSJ:13}.

\begin{theorem}\label{thm:QP_optimal_preconditioning}
Let $R_q R_q^\top = Q^{-1}$ be the Choleski factorization of $Q^{-1}$ and $P\in \R{n \times n-s}$ be a matrix whose columns are orthonormal vectors spanning $\Range{R_q^\top A^\top}$ with $s$ being the dimension of $\Null{A}$ and let $\lambda_n(LAQ^{-1}A^\top L)$ and $\lambda_1(LAQ^{-1}A^\top L)$ be the largest and smallest nonzero eigenvalues of $LAQ^{-1}A^\top L$. The diagonal scaling matrix $L^\star \in\PD{m}$ that minimizes the eigenvalue ratio $\lambda_n(LAQ^{-1}A^\top L)/\lambda_1(LAQ^{-1}A^\top L)$ can be obtained by solving the convex problem
\begin{equation}\label{eqn:QP_optimal_scaling_convex}
\begin{aligned}
\begin{array}{ll}
%
\underset{{t\in\R{},\;w\in\R{m}}}{\mbox{minimize}} & t\\
\mbox{subject to} & W=\mbox{diag}(w),\; w>0,\\
& tI -  R_q^\top A^\top W A R_q  \in\PSD{n},\\
&  P^\top (R_q^\top A^\top W A R_q -  I )P  \in\PSD{n-s},
\end{array}
\end{aligned}
\end{equation}
and setting $L^\star = W^{\star^{1/2}}$.

\end{theorem}


So far, we characterized the convergence factor of the ADMM algorithm based on general properties of the sequence $\{ F^k v^k\}$. However, if we a priori know which constraints will be active during the ADMM iterations, our parameter selection rules~\eqref{eqn:QP_optimal_factor} and~\eqref{eqn:QP_relaxation_optimal_factor} may not be optimal. To illustrate this fact, we will now analyze the two extreme situations where no and all constraints are active in each iteration and derive the associated optimal ADMM parameters.
\subsection{Special cases of quadratic programming}
The first result deals with the case where the constraints of~\eqref{eqn:Quadratic_problem} are never active. This could happen, for example, if we use the constraints to impose upper and lower bounds on the decision variables, and use very loose bounds.

\begin{proposition}\label{prop:QP_when_careful_1}
Assume that ${F^{k+1} = F^{k} = -I}$ for all epochs $k\in \R{}_{+}$ in~\eqref{eqn:Quadratic_admm_iterations_reformed} and~\eqref{eqn:Quadratic_admm_iterations_reformed_relaxation}.  Then the modified ADMM algorithm~\eqref{eqn:QP_Fv_sequence_relaxation} attains its minimal convergence factor for the parameters
\begin{align}
\alpha = 1,\quad \rho \rightarrow 0.
\end{align}
In this case \eqref{eqn:QP_Fv_sequence_relaxation} coincide with~\eqref{eqn:QP_Fv_sequence} and their convergence factor is minimized: $\zeta = \zeta_R \rightarrow 0$.
\end{proposition}
The next proposition addresses another extreme scenario when the ADMM iterates are operating on the active set of the quadratic program~\eqref{eqn:Quadratic_problem}.
\begin{proposition}\label{prop:QP_when_careful_2}
Suppose that $F^{k+1} = F^{k} = I$ for all $k\in \R{}_+$~in~\eqref{eqn:Quadratic_admm_iterations_reformed} and~\eqref{eqn:Quadratic_admm_iterations_reformed_relaxation}. Then the relaxed ADMM algorithm~\eqref{eqn:QP_Fv_sequence_relaxation} attains its minimal convergence factor for the parameters
\begin{align}
\alpha = 1,\quad \rho \rightarrow \infty.
\end{align}
In this case \eqref{eqn:QP_Fv_sequence_relaxation} coincides with~\eqref{eqn:QP_Fv_sequence} and their convergence factors are minimized: $\zeta = \zeta_R \rightarrow 0$.
\end{proposition}

It is worthwhile to mention that when~\eqref{eqn:Quadratic_problem} is defined so that its constraints are active (inactive) then the $s^{k}$ ($r^k$) residuals of the ADMM algorithm remain zero for all  $k \geq 2$ updates.

\section{Numerical examples}
\label{sec:qp_evaluation}
In this section, we  evaluate our parameter selection rules on numerical examples. First, we illustrate the convergence factor of ADMM and gradient algorithms for a family of $\ell_2$-regularized quadratic problems. These examples demonstrate that the ADMM method converges faster than the gradient method for certain ranges of the regularization parameter $\delta$, and slower for other values. Then, we consider QP-problems and
compare the performance of the over-relaxed ADMM algorithm with an alternative accelerated ADMM method presented in~\cite{GOS:2012}. The two algorithms are also applied to a Model Predictive Control (MPC) benchmark where QP-problems are solved repeatedly over time for fixed matrices $Q$ and $A$ but varying vectors $q$ and $b$.

\subsection{$\ell_2$-regularized quadratic minimization via ADMM}
\begin{figure}[h]
\centering
  \includegraphics[width=.4\columnwidth]{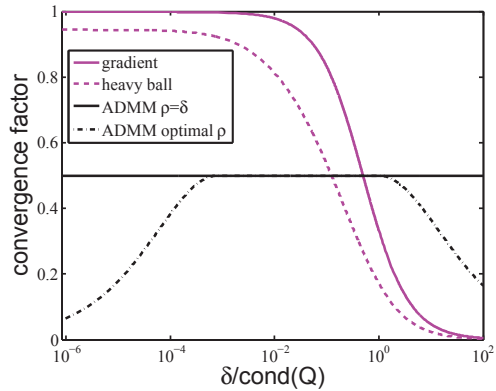}
\caption{Convergence factor of the ADMM, gradient, and heavy-ball methods for $\ell_2$ regularized minimization with fixed $Q$-matrix and different values of the regularization parameter $\delta$.}
\label{fig:l2_rate}
\end{figure}

We consider $\ell_2$-regularized quadratic minimization problem~\eqref{cor:L2:standard} for a $Q\in \mathcal{S}_{++}^{100}$ with condition number $1.2 \times 10^3$ and for a range of regularization parameters $\delta$.
Fig.~\ref{fig:l2_rate} shows  how the optimal convergence factor of ADMM depends on $\delta$. The results are shown for
 two step-size rules: $\rho=\delta$ and $\rho=\rho^\star$ given in~\eqref{eqn:ADMM_L2_optimal_step-size}. For comparison, the gray and dashed-gray curves show the optimal convergence factor of the gradient method
\begin{align*}
&x^{k+1}=x^k - \gamma (Qx^k + q + \delta x^k),\\
\intertext{with step-size $\gamma<2/(\lambda_n(Q)+\delta)$ and  a multi-step gradient iterations on the form}
&x^{k+1}=x^k - a (Q x^k + q + \delta x^k) + b (x^k- x^{k-1}),
\end{align*}
This latter algorithm is known as the heavy-ball method and significantly outperforms the standard gradient method on ill-conditioned problems~\cite{polyak}. The algorithm has two parameters: $a<2(1+b)/(\lambda_n(Q)+\delta)$, and $b\in [0,1]$. For our problem, since the cost function is quadratic and its Hessian $\nabla^2 f(x)= Q + \delta I$ is bounded between $l=\lambda_1(Q)+\delta$ and $u=\lambda_n(Q)+\delta$, the optimal step-size for the gradient method is $\gamma^{\star}=2/(l+u)$ and the optimal parameters for the heavy-ball method are $a^\star=4/(\sqrt{l}+\sqrt{u})^2$, and $b^\star = (\sqrt{u}-\sqrt{l})^2/(\sqrt{l}+\sqrt{u})^2$\cite{polyak}.

Figure~\ref{fig:l2_rate} illustrates the convergence properties of the ADMM method under both step-size rules. The optimal step-size rule gives significant speedups of the ADMM for small or large values of the regularization parameter $\delta$. This phenomena can be intuitively explained based on the interplay of the two parts of the objective function in~\eqref{eqn:L2_formulation}. For extremely small values of $\delta$, one sees that the $x$-th part of the objective is becoming dominant compared to $z$-th part. Consequently, using the optimal step-size in~\eqref{eqn:ADMM_L2_optimal_step-size}, $z$- is dictated to quickly follow the value of $x$-update. A similar reasoning holds when $\delta$ is large, in which the $x$- has to obey the $z$-update.

It is interesting to observe that ADMM outperforms the gradient and heavy-ball methods for small $\delta$ (an ill-conditioned problem), but actually performs worse as $\delta$ grows large (\ie when the regularization makes the overall problem well-conditioned). It is noteworthy that the relaxed ADMM method solves the same problem in one step (convergence factor $\zeta^\star_R=0$).
\subsection{Quadratic programming via ADMM}
Next, we evaluate our step-size rules for ADMM-based quadratic programming and compare their performance with that of other accelerated ADMM variants from the literature.

\subsubsection{Accelerated ADMM}
One recent proposal for accelerating the ADMM-iterations is called \emph{fast-ADMM}~\cite{GOS:2012} and consists of the following iterations
\begin{align}
\label{eqn:admm_nesterov_acceleration}
\begin{array}[c]{ll}
	x^{k+1} &= \underset{x}{\operatorname{argmin}}\, L_{\rho}(x,\hat{z}^{k}, \hat{u}^{k}), \\
	z^{k+1} &= \underset{z}{\operatorname{argmin}}\, L_{\rho}(x^{k+1}, z, \hat{u}^k), \\
	u^{k+1} &= \hat{u}^{k} + Ax^{k+1}+Bz^{k+1}-c, \\
    \hat{z}^{k+1} &= \alpha^k z^{k+1}+(1-\alpha^k)z^k,\\
    \hat{u}^{k+1} &= \alpha^k u^{k+1} + (1-\alpha^k)u^k.
    \end{array}
\end{align}
 The relaxation parameter $\alpha^k$ in the fast-ADMM method is defined based on the Nesterov's order-optimal method~\cite{Nesterov03} combined with an innovative restart rule where $\alpha^k$ is given by
  \begin{align}\label{eqn:fast_ADMM_restart_rule}
  \alpha^k =\left\{
  \begin{array}[c]{ll}
  1+\dfrac{\beta^k -1}{\beta^{k+1}} & \operatorname{if}\, \dfrac{\max(\Vert r^k \Vert, \Vert s^k\Vert)}{\max(\Vert r^{k-1}\Vert, \Vert s^{k-1}\Vert)}< 1,  \\
  1 &\mbox{otherwise},
  \end{array}\right.
  \end{align}
  where $\beta^1=1$, and $\beta^{k+1}= \dfrac{1+\sqrt{1+4{\beta^k}^2}}{2}$ for $k>1$.
  The restart rule assures that~\eqref{eqn:admm_nesterov_acceleration} is updated in the descent direction with respect to the primal-dual residuals.

\begin{figure}[h]
        \begin{center}
	\subfigure[$n=100,\quad m=50$.]{\label{fig:qp_feasible}
        \includegraphics[width=.4\hsize]{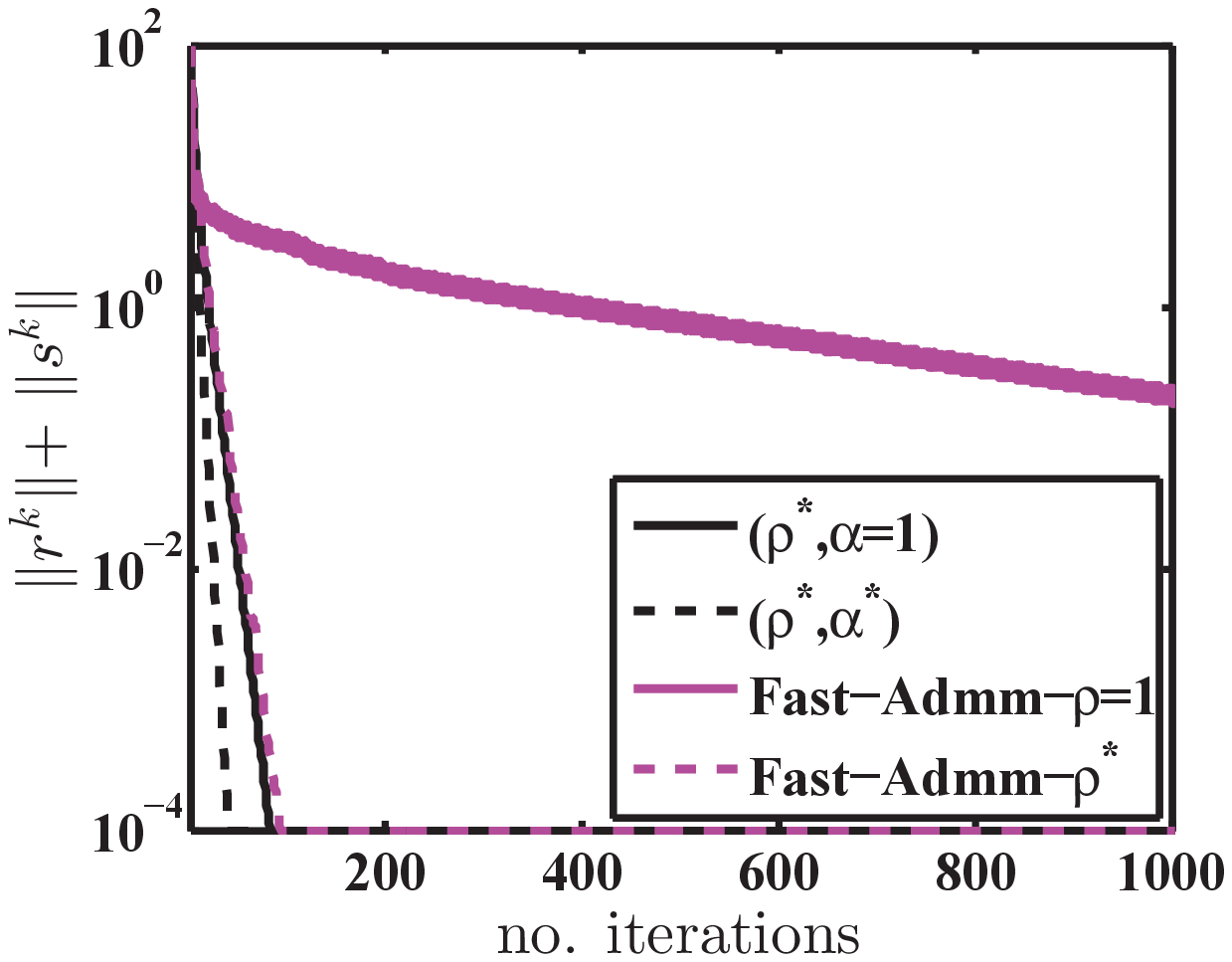}}
 	\subfigure[$n=100,\quad m=200$]{
        \includegraphics[width=.4\hsize]{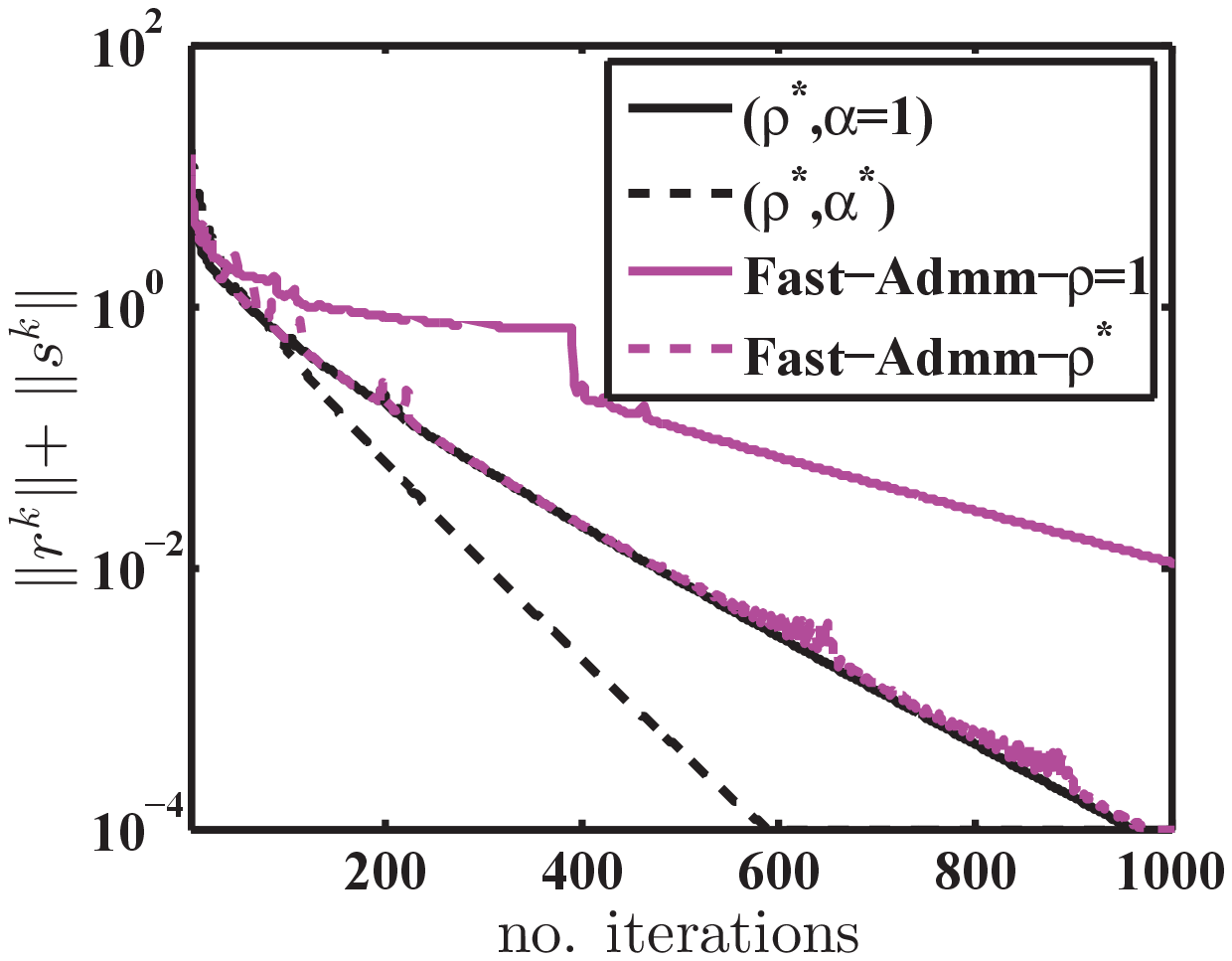}
        \label{fig:qp_active}}
\end{center}
    \caption{Convergence of primal plus dual residuals of four ADMM algorithms with $n$ decision variables and $m$ inequality constraints.}
\label{fig:qp_comparison}
\end{figure}

To compare the performance of the over-relaxed ADMM iterations with our proposed parameters to that of fast-ADMM, we conducted several numerical examples.
For the first numerical comparison, we generated several instances of~\eqref{eqn:Quadratic_problem}; Figure~\ref{fig:qp_comparison} shows the results for the two representative examples. In the first case, $A\in \mathcal{R}^{50\times 100}$ and $Q\in \mathcal{S}_{++}^{100}$ with condition number $1.95 \times 10^{3}$; $32$ constraints are active at the optimal solution. In the second case,
$A\in \mathcal{R}^{200\times 100}$ and $Q\in \mathcal{S}_{++}^{100}$, where  the condition number of $Q$ is $7.1\times 10^3$. The polyhedral constraints correspond to random box-constraints, of which $66$ are active at optimality. We evaluate for four algorithms: the ADMM iterates in~\eqref{eqn:Quadratic_admm_iterations_relaxation} with and without over-relaxation and the corresponding tuning rules developed in this paper, and the fast-ADMM iterates~\eqref{eqn:admm_nesterov_acceleration} with $\rho=1$ as proposed by~\cite{GOS:2012} and $\rho=\rho^\star$ of our paper. The convergence of corresponding algorithms in terms of the summation of primal and dual residuals $\Vert r^k\Vert+ \Vert s^k\Vert$ are depicted in Fig.~\ref{fig:qp_comparison}. The plots exhibit a significant improvement of our tuning rules compared to the fast-ADMM algorithm.

To the best of our knowledge, there are currently no results about optimal step-size parameters for the fast-ADMM method. However, based on our numerical investigations, we observed that the performance of fast-ADMM algorithm significantly improved by employing our optimal step-size $\rho^\star$ (as illustrated in~\ref{fig:qp_comparison}). In the next section we perform another comparison between three algorithms, using the optimal $\rho$-value for fast-ADMM obtained by an extensive search.

%
%
%
%

\begin{figure}[h]
\begin{center}
    \subfigure[$\alpha = 1$, $L=I$]{\label{fig:QP_MPCa}
        \includegraphics[width=.4\hsize ]{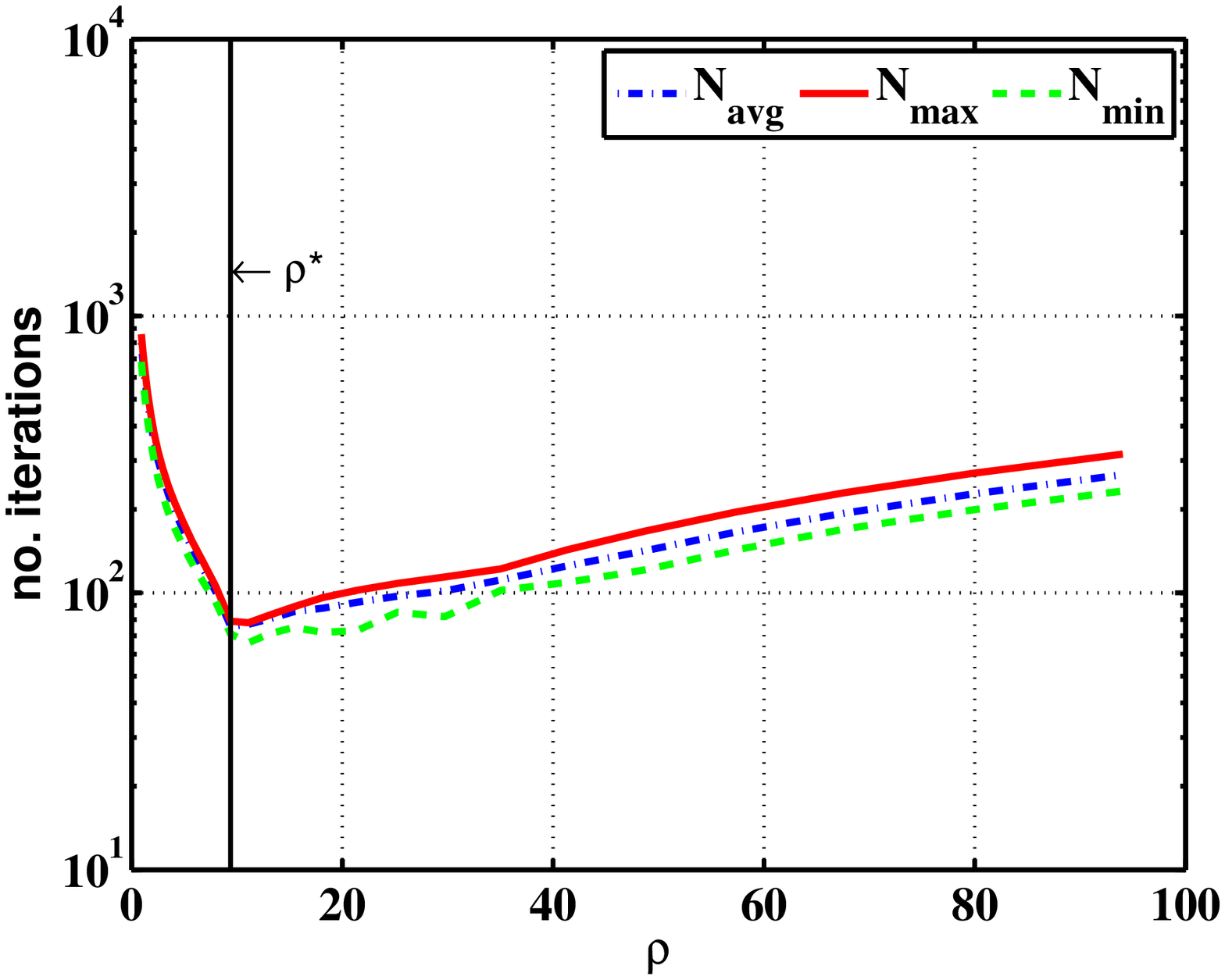}}%
    \subfigure[$\alpha = 2$, $L=I$]{\label{fig:QP_MPCb}
        \includegraphics[width=.4\hsize]{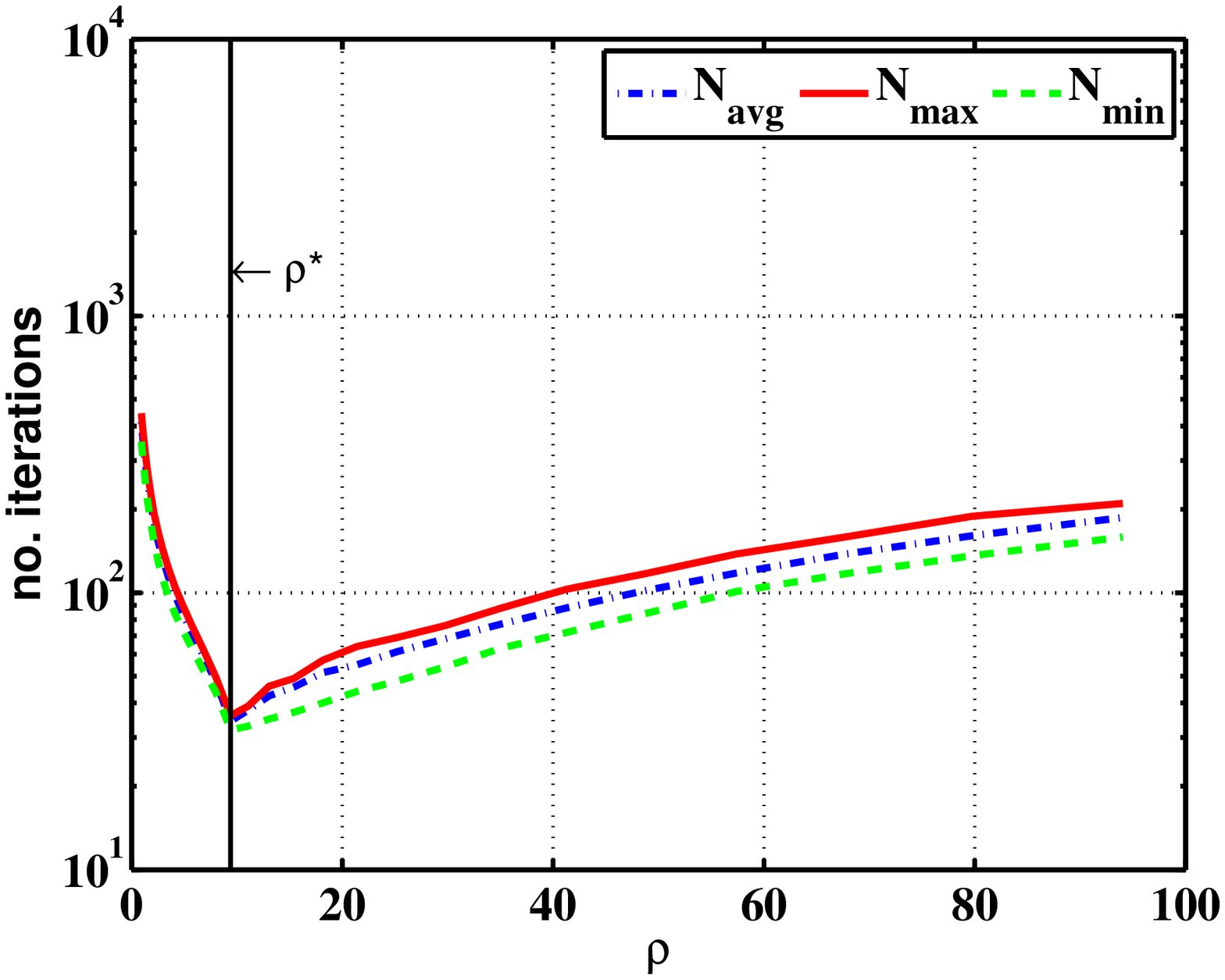}}
        \subfigure[$\alpha = 1$, $L=L^\star$]{
        \includegraphics[width=.4\hsize]{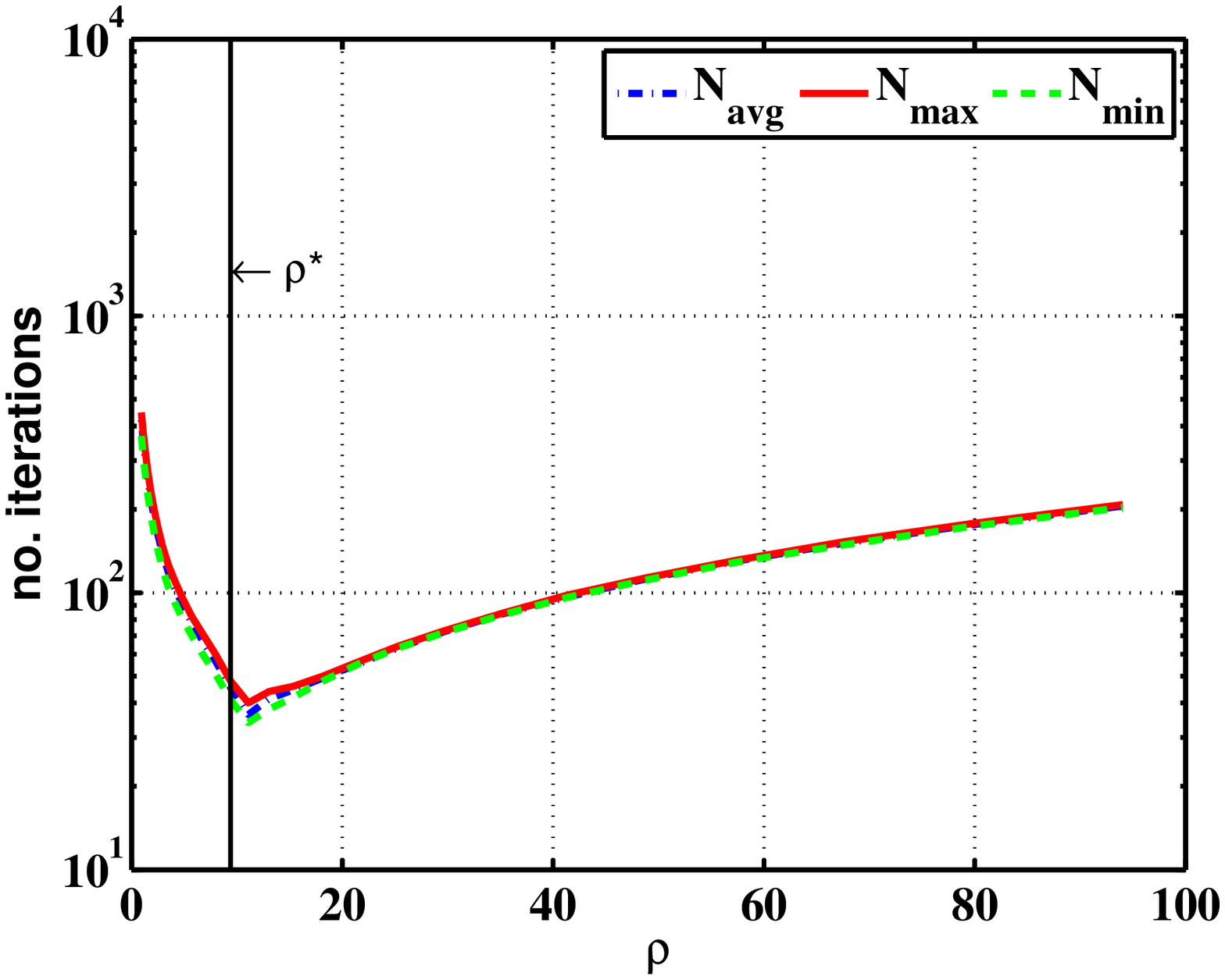}
        \label{fig:QP_MPCc}}%
    \subfigure[$\alpha = 2$, $L=L^\star$]{
        \includegraphics[width=.4\hsize]{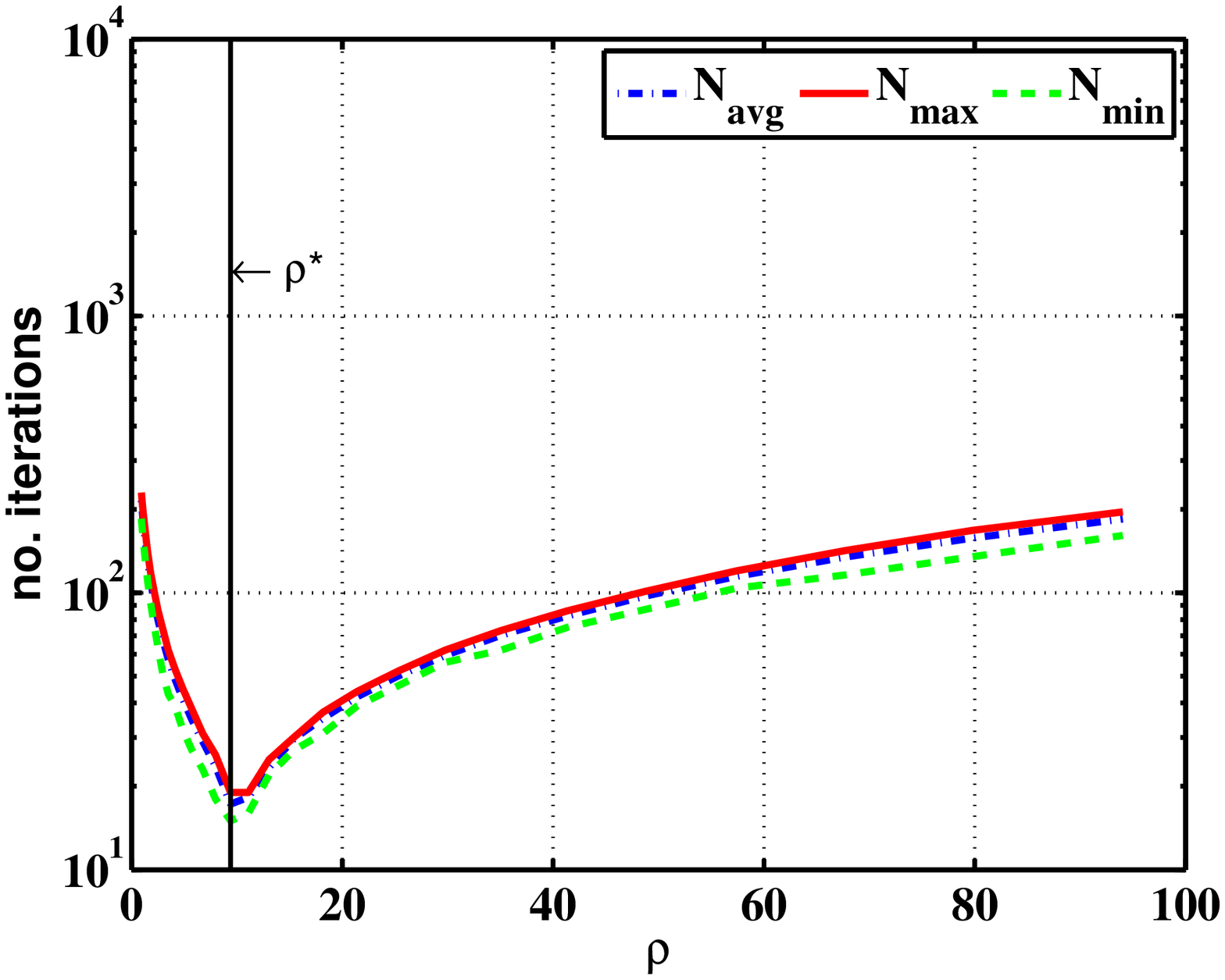}
        \label{fig:QP_MPCd}}
\end{center}
    \caption{Number of iterations $k:\,\max\{\|r^k\|,\, \|s^k\|\} \leq 10^{-5}$ for ADMM applied to the MPC problem for different initial states $x_0$. The dashed green line  denotes the minimum number of iterations taken over all the initial states, the dot-dashed blue line corresponds to the average, while the red solid line represents the maximum number of iterations.}\label{fig:QP_MPC}
\end{figure}

\subsubsection{Model Predictive Control}
Consider the discrete-time linear system
\begin{align}\label{eq:control_system}
x_{t+1} &= Hx_t + J u_t + J_r r,
\end{align}
where $t\geq 0$ is the time index, $x_t \in \R{n_x}$ is the state, $u_t \in \R{n_u}$ is the control input,  $r \in \R{n_r}$ is a constant reference signal, and $H\in \R{n_x \times n_x}$, $J\in \R{n_x \times n_u}$, and $J_r\in \R{n_x \times n_r}$ are fixed matrices.
Model predictive control aims at solving the following optimization problem
\begin{equation}
\begin{aligned}
\label{eqn:MPC_problem_1}
	\begin{array}[c]{ll}
	\underset{\{u_i\}_0^{N_p-1}}{\mbox{minimize}} &  \dfrac{1}{2}\sum_{i=0}^{N_p - 1}(x_i-x_r)^{\top}Q_x (x_i-x_r) + (u_i-u_r)^{\top} R (u_i-u_r) + (x_{N_p}-x_r)^\top Q_N (x_{N_p}-x_r)\\
	\mbox{subject to}& x_{t+1} = Hx_t + J u_t + J_r r\quad \forall t,\\
                     & x_t\in\mathcal{C}_x\quad \forall t,\\
                     & u_t\in\mathcal{C}_u\quad \forall t,
	\end{array}
\end{aligned}
\end{equation}
where $x_0$, $x_r$, and $u_r$ are given, $Q_x\in \PD{n_x}$, $R\in \PD{n_u}$, and $Q_N\in \PD{n_x}$ are the state, input, and terminal costs, and the sets $\mathcal{C}_x$ and $\mathcal{C}_u$ are convex. Suppose that the sets $\mathcal{C}_x$ and $\mathcal{C}_u$ correspond to component-wise lower and upper bounds, i.e., $\mathcal{C}_x = \{x\in\R{n_x} \vert 1_{n_x}\bar{x}_{min}\leq x \leq 1_{n_x}\bar{x}_{max}\}$ and $\mathcal{C}_u = \{u\in\R{n_u} \vert 1_{n_u}\bar{u}_{min}\leq u \leq 1_{n_u}\bar{u}_{max}\}$.
Defining $\chi = [x_1^\top \, \dots \, x_{N_p}^\top]^\top$, $\upsilon=[u_0^\top \, \dots \, u_{N_p-1}^\top]^\top$, $\upsilon_r=[r^\top \, \dots \, r^\top]^\top$,~\eqref{eq:control_system} can be rewritten as $\chi = \Theta x_0 + \Phi\upsilon + \Phi_r\upsilon_r$. The latter relationship can be used to replace $x_t$ for $t=1,\dots,N_p$ in the optimization problem, yielding the following QP:
\begin{equation}
\begin{aligned}
\label{eqn:MPC_problem_4}
	\begin{array}[c]{ll}
	\underset{\upsilon}{\mbox{minimize}} &  \dfrac{1}{2}\upsilon^\top Q \upsilon + q^\top \upsilon\\
	\mbox{subject to} & A\upsilon \leq b ,
	\end{array}
\end{aligned}
\end{equation}
where
\begin{equation}
\bar{Q} =
\begin{bmatrix}
I_{N_p-1} \otimes Q_x & 0\\
0   &             Q_N
\end{bmatrix},\quad
\bar{R} = I_{N_p} \otimes R,\quad
\begin{aligned}
A =
\begin{bmatrix}
\Phi\\
-\Phi\\
I\\
-I
\end{bmatrix},\quad
b=\begin{bmatrix}
1_{n_x N_p} \bar{x}_{max} - \Theta x_0-\Phi_r\upsilon_r\\
1_{n_x N_p} \bar{x}_{min} + \Theta x_0+\Phi_r\upsilon_r\\
1_{n_u N_p} \bar{u}_{max}\\
1_{n_u N_p} \bar{u}_{min}
\end{bmatrix},
\end{aligned}
\end{equation}
and $Q=\bar{R} + \Phi^\top\bar{Q}\Phi$ and $q^\top = x_0^\top\Theta^\top \bar{Q}\Phi + \upsilon_r^\top \Phi_r^\top \bar{Q}\Phi - x_r^\top\left( 1_{N_p}^\top\otimes I_{n_x} \right)\bar{Q}\Phi - u_r^\top \left( 1_{N_p}^\top\otimes I_{n_u} \right)\bar{R}$.

Below we illustrate the MPC problem for the quadruple-tank process~\cite{Johansson2000}. The state of the process $x\in\R{4}$ corresponds to the water levels of all tanks, measured in centimeters. The plant model was linearized at a given operating point and discretized with a sampling period of $2\,s$. The MPC prediction horizon was chosen as $N_p = 5$. A constant reference signal was used, while the initial condition $x_0$ was varied to obtain a set of MPC problems with different non-empty feasible sets and linear cost terms. In particular, we considered initial states of the form $x_0 = [x_1\, x_2\, x_3\, x_4]^\top$ where $x_i \in \{10,\; 11.25,\;  12.5,\;   13.75,\;   15\}$ for $i=1,\dots,4$. Out of the possible $625$ initial values, $170$ yields feasible QPs (each with $n=10$ decision variables and $m=40$ inequality constraints). We have made these QPs publically available as a {\tt MATLAB} formatted binary file~\cite{mpc_dataset}. To prevent possible ill-conditioned QP-problems, the constraint matrix $A$ and vector $b$ were scaled so that each row of $A$ has unit-norm.
 %

Fig.~\ref{fig:QP_MPC} illustrates the convergence of the ADMM iterations for the $170$ QPs as a function of the step-size $\rho$, scaling matrix $L$, and over-relaxation factor $\alpha$. Since $A^\top$ has a non-empty null-space, the step-size $\rho^\star$ was chosen heuristically based on Theorem~\ref{thm:QP_optimal_factor} as $\rho^\star = 1/\sqrt{\lambda_{1}(AQ^{-1}A^\top)\lambda_{n}(AQ^{-1}A^\top)}$, where $\lambda_1(AQ^{-1}A^\top)$ is the smallest nonzero eigenvalue of $AQ^{-1}A^\top$.
As shown in Fig.~\ref{fig:QP_MPC}, our heuristic step-size $\rho^\star$ results in a number of iterations close to the empirical minimum. Moreover, performance is improved by choosing $L=L^\star$ and $\alpha=2$.

\begin{figure}[tb]
\centering
\includegraphics[width=.4\hsize]{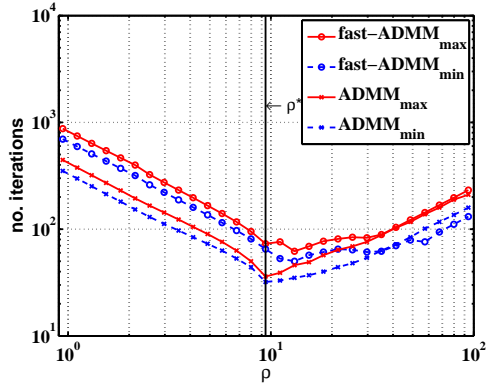}
    \caption{Number of iterations $k:\,\max\{\|r^k\|,\, \|s^k\|\} \leq 10^{-5}$ for ADMM with $L=I$ and $\alpha=2$ and fast-ADMM algorithms applied to the MPC problem for different initial states $x_0$. The line in blue denotes the minimum number of iterations taken over all the initial states, while the red line represents the maximum number of iterations.}\label{fig:QP_MPC_FADMMvsADMM}
\end{figure}

The performance of the Fast-ADMM and ADMM algorithms is compared in Fig.~\ref{fig:QP_MPC_FADMMvsADMM} for $L=I$ and $\alpha=2$. The ADMM algorithm with the optimal over-relaxation factor $\alpha=2$ uniformly outperforms the Fast-ADMM algorithm, even with suboptimal scaling matrix $L$.

\subsubsection{Local convergence factor}
To illustrate our results on the slow local convergence of ADMM, we consider a QP problem of the form~\eqref{eqn:MPC_problem_4} with
\begin{equation}
\label{eqn:QP:slow_covergence_example}
\begin{aligned}
Q&=
\begin{bmatrix}
40.513  &  0.069\\
0.069 &  40.389
\end{bmatrix},\quad q=0\\
A&=\begin{bmatrix}
-1   &      0 \\
0   &-1         \\
0.1151   & 0.9934\\
\end{bmatrix},\quad
b=
\begin{bmatrix}
6\\
6\\
-0.3422
\end{bmatrix}.
\end{aligned}
\end{equation}
%
\begin{figure}[tb]
\begin{center}
    \subfigure[a]{\includegraphics[width=.4\hsize]{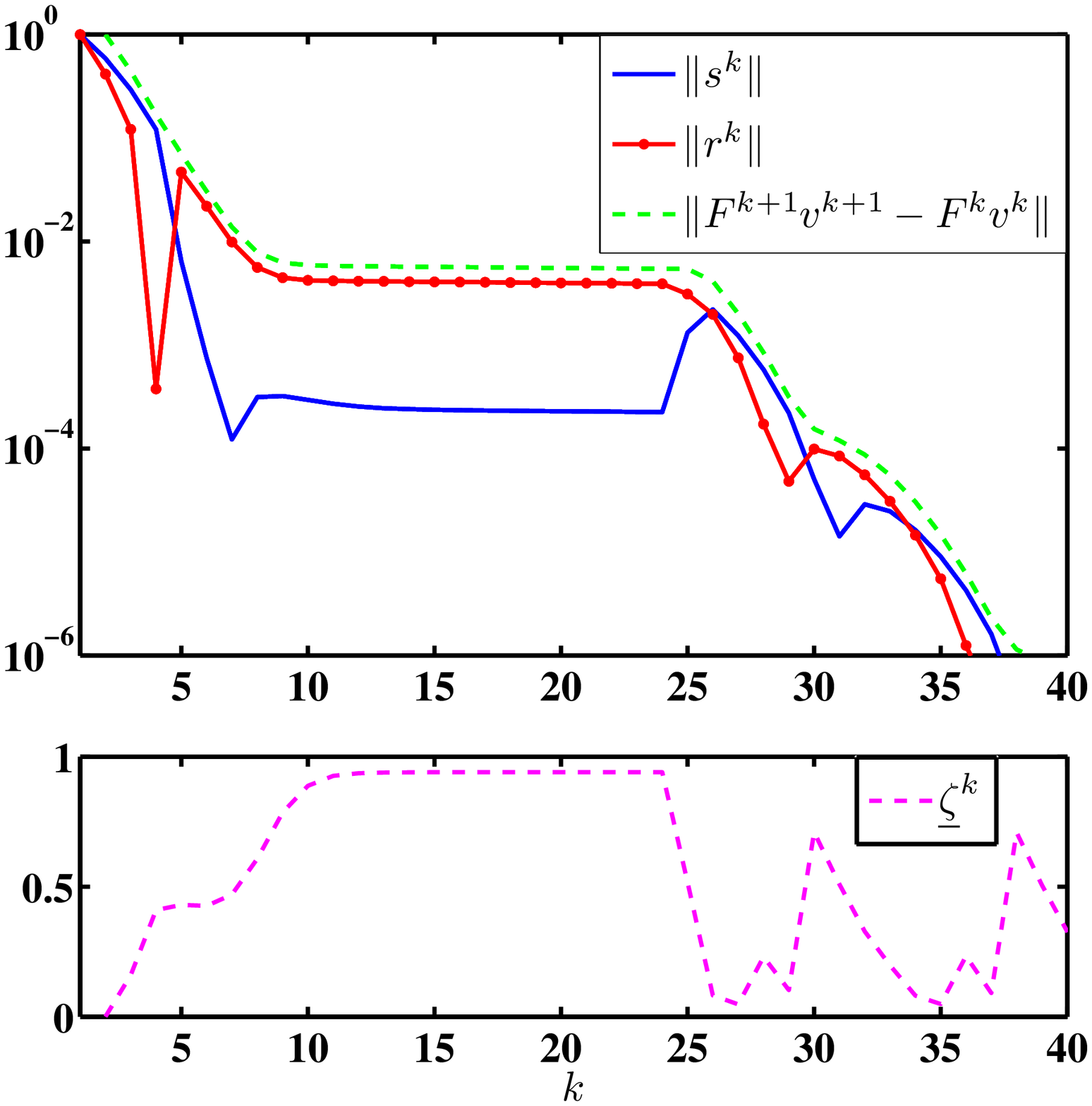}
        \label{fig:QP_slow_iters}}
    \subfigure[b]{
        \includegraphics[width=.4\hsize]{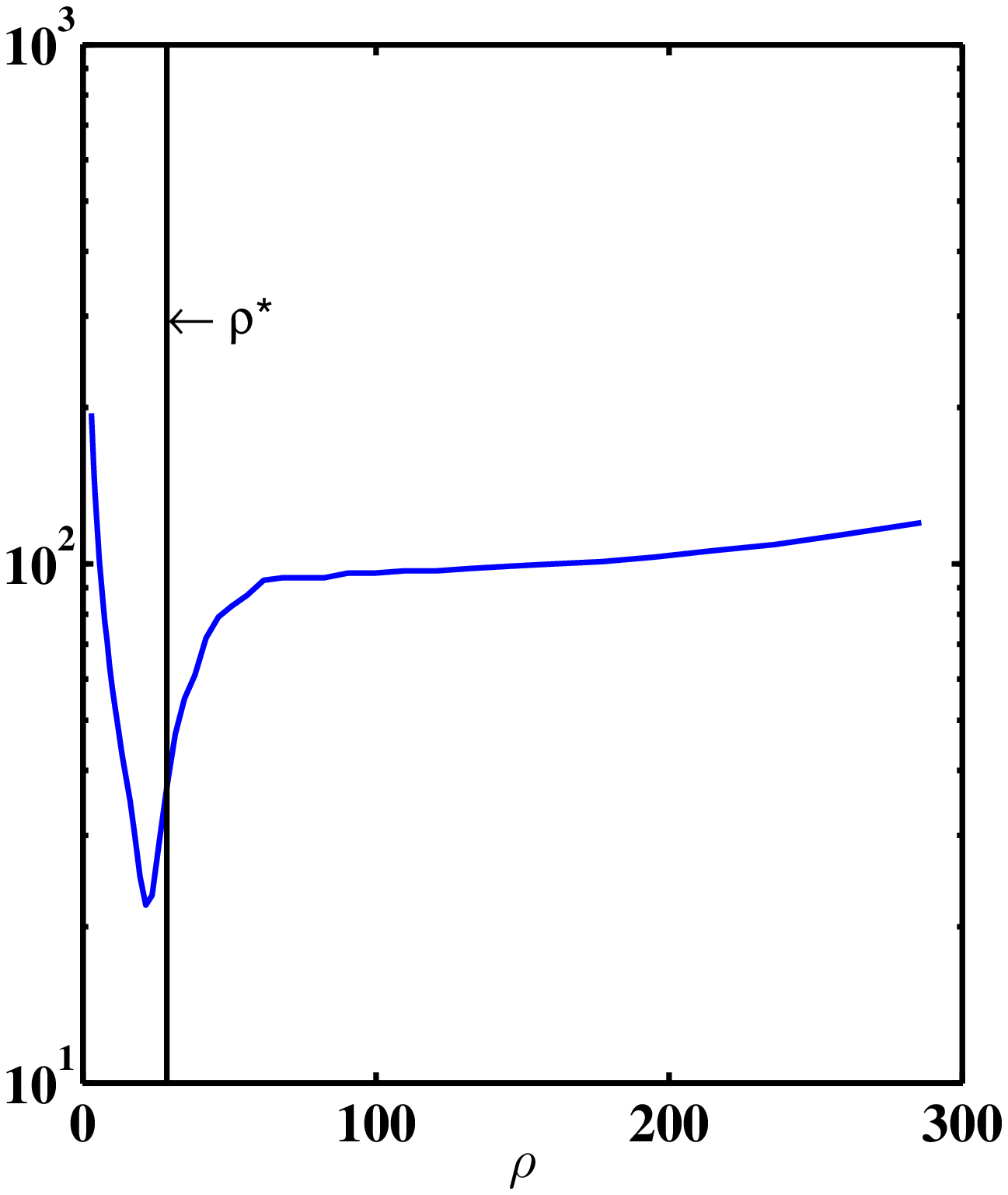}}
\end{center}
    \caption{Slow convergence of ADMM algorithm for the example in~\eqref{eqn:QP:slow_covergence_example} with $\alpha=1$ and $L=I$. The residuals $r^k$, $s^k$, and $F^{k+1}v^{k+1} - F^kv^k$ and the lower bound on the convergence factor $\underline\zeta^k$ are shown in the left, while the number of iterations for $\rho\in[0.1\rho^\star \; 10\rho^\star]$ are shown in the right.
     \label{fig:QP_slow}}
\end{figure}

The ADMM algorithm was applied to the former optimization problem with $\alpha = 1$ and $L=I$. Given that the nullity of $A$ is not $0$, the step-size was chosen heuristically based on Theorem~\ref{thm:QP_optimal_factor} as $\rho^\star = 1/\sqrt{\lambda_{1}(AQ^{-1}A^\top)\lambda_{n}(AQ^{-1}A^\top)} = 28.6$ with $\lambda_1(AQ^{-1}A^\top)$ taken to be the smallest nonzero eigenvalue of $AQ^{-1}A^\top$. The resulting residuals are shown in Fig.~\ref{fig:QP_slow}, together with the lower bound on the convergence factor $\underline\zeta$ evaluated at each time-step. As expected from the results in Theorem~\ref{thm:QP:linear_rate}, the residual $F^{k+1}v^{k+1} - F^{k}v^{k}$ is monotonically decreasing. However, as illustrated by $\underline\zeta^k$, the lower bound on the convergence factor from Theorem~\ref{thm:QP_slow}, the auxiliary residual $F^{k+1}v^{k+1} - F^{k}v^{k}$ and the primal-dual residuals show a convergence factor close to $1$ over several time-steps. The heuristic step-size rule performs reasonably well as illustrated in the right subplot of Fig.~\ref{fig:QP_slow}.

\section{Conclusions and Future Work}
We have studied optimal parameter selection for the alternating direction method of multipliers for two classes of quadratic problems: $\ell_2$-regularized quadratic minimization and quadratic programming under linear inequality constraints. For both problem classes, we established global convergence of the algorithm at linear rate and provided explicit expressions for the parameters that ensure the smallest possible convergence factors. We also considered iterations accelerated by over-relaxation, characterized the values of the relaxation parameter for which the over-relaxed iterates are guaranteed to improve the convergence times compared to the non-relaxed iterations, and derived jointly optimal step-size and relaxation parameters. We validated the analytical results on numerical examples and demonstrated superior performance of the tuned ADMM algorithms compared to existing methods from the literature. As future work, we plan to extend the analytical results for more general classes of objective functions.

\section*{Acknowledgment}
\addcontentsline{toc}{section}{Acknowledgment}                             
The authors would like to thank Pontus Giselsson, Themistoklis Charalambous, Jie Lu, and Chathuranga Weeraddana
for their valuable comments and suggestions to this manuscript.

\bibliographystyle{IEEEtran}
\bibliography{admmbib}

\appendix
\section{Proofs}    

\subsection{Proof of Theorem~\ref{thm:L2:standard}}
From Proposition~\ref{prop:1}, the variables $x^k$ and $z^k$ in iterations (\ref{eqn:ADMM_L2_iterations}) converge to the optimal values $x^\star$ and $z^\star$ of~\eqref{eqn:L2_formulation} if and only if the spectral radius of the matrix $E$ in  (\ref{eqn:ADMM_L2_matrix}) is less than one. To express the eigenvalues of $E$ in terms of the eigenvalues of $Q$, let $\lambda_i(Q), i=1, \dots, n$ be the eigenvalues of $Q$ sorted in ascending order. Then, the eigenvalues $\zeta(\rho, \lambda_i(Q))$ of $E$ satisfy
\begin{align}
\zeta(\rho, \lambda_i(Q))
=\dfrac{\rho^2+ \lambda_i(Q) \delta}{\rho^2+\lambda_i(Q) \delta+(\lambda_i(Q)+\delta)\rho}.\label{eqn:fexpressions}
\end{align}
Since $\lambda_i(Q), \rho, \delta \in \mathcal{R}_{++}$, we have $0 \leq \zeta(\rho, \lambda_i(Q)) <1$ for all $i$, which ensures convergence.

To find the optimal step-size parameter and the associated convergence factor $(\rho^{\star}, \zeta^\star)$, note that, for a fixed $\rho$, the convergence factor $\zeta(\rho) = \max_{e^k} \|e^{k+1}\|/\|e^k\|$ corresponds to the spectral radius of $E$, \ie  $\zeta(\rho)=\max_i\left\{ \zeta(\rho,\lambda_i(Q))\right\}$. It follows that the optimal pair $(\rho^{\star}, \zeta^\star)$ is given by
\begin{equation}
\begin{aligned}
\label{eqn:ADMM_L2_optimal_step-size_optProblem}
\rho^{\star} = \underset{\rho}{\mbox{argmin }} \max_i\left\{ \zeta(\rho,\lambda_i(Q))\right\},\quad
\zeta^{\star} = \max_i\left\{ \zeta(\rho^\star,\lambda_i(Q))\right\}.
\end{aligned}
\end{equation}

From ~(\ref{eqn:fexpressions}), we can see that $\zeta(\rho,\lambda_i(Q))$ is monotone decreasing in $\lambda_i(Q)$ when $\rho>\delta$ and monotone increasing when $\rho<\delta$. Hence, we consider these two cases separately.

When $\rho>\delta$, the largest eigenvalue of $E$ is given by $\zeta(\rho, \lambda_1(Q))$ and $\rho^{\star} = \mbox{argmin}_{\rho}\zeta(\rho,\lambda_1(Q))$.
By the first-order optimality conditions and the explicit expressions in~\eqref{eqn:fexpressions} we have
\begin{align*}
	\rho^{\star}&=\sqrt{\delta \lambda_1(Q)}, \quad \zeta^\star = \zeta(\rho^{\star}, \lambda_1(Q)) = (1+\dfrac{\delta+\lambda_1(Q)}{2\sqrt{\delta \lambda_1(Q)}})^{-1}.
\end{align*}
However, this value of $\rho$ is larger than $\delta$ only if $\delta<\lambda_1(Q)$.  When $\delta\geq \lambda_1(Q)$, the assumption that $\rho>\delta$ implies that $0\leq (\rho-\delta)^2 \leq (\rho-\delta)(\rho-\lambda_1(Q))$, so
\begin{align*}
	&\zeta(\rho,\lambda_1(Q)) =\dfrac{\rho^2+ \lambda_i(Q) \delta}{\rho^2+\lambda_i(Q) \delta+(\lambda_i(Q)+\delta)\rho}\geq\\
& \dfrac{\rho^2+\lambda_1(Q)\delta}{\rho^2 +\lambda_1(Q)\delta +(\lambda_1(Q)+\delta)\rho + (\rho-\delta)(\rho-\lambda_1(Q))}=\dfrac{1}{2}.
\end{align*}
Since $\rho=\delta$ attains $\zeta(\delta,\lambda_1(Q))=1/2$ it is optimal.

A similar argument applies to $\rho<\delta$. In this case, $\max_i \zeta(\rho,\lambda_i(Q)) = \zeta(\rho, \lambda_n(Q))$ and
 when $\delta>\lambda_n(Q)$, $\rho^{\star}=\sqrt{\delta \lambda_n(Q)}$ is the optimal step-size  and the associated convergence factor is \[\zeta^\star =  \left(1+\dfrac{\delta+\lambda_n(Q)}{2\sqrt{\delta \lambda_n(Q)}}\right)^{-1}.\] For $\delta\leq \lambda_n(Q)$, the requirement that $\rho<\delta$ implies the inequalities $0\leq (\delta-\rho)^2 \leq (\lambda_n(Q)-\rho)(\delta-\rho)$ and that $\zeta(\rho,\lambda_n(Q)) \geq \dfrac{1}{2}$,
which leads to $\rho=\delta$ being optimal. 

\subsection{Proof of Corollary~\ref{cor:L2:standard}}
The proof is a direct consequence of evaluating \eqref{eqn:fexpressions} at $\rho=\delta$ for $i=1,\dots,n$.
\subsection{Proof of Theorem~\ref{thm:L2:relaxed}}
The $z$-update in~\eqref{eqn:ADMM_L2_iterations_relaxation} implies that $\mu^k = (\delta+\rho)z^{k+1}- \rho (\alpha x^{k+1}+ (1-\alpha)z^k)$, and that the $\mu$-update in~\eqref{eqn:ADMM_L2_iterations_relaxation} can be written as $\mu^{k+1} = \delta z^{k+1}$.
Similarly to the analysis of the previous section, inserting the $x$-update into the $z$-update, we find
\begin{align*}
%
z^{k+1} = \underset{E_R}{\underbrace{\dfrac{1}{\delta+\rho}\left(\delta I + \rho \left(\alpha(\rho-\delta)\left(Q+\rho I\right)^{-1}+ (1-\alpha)I\right)\right)} }z^k - \dfrac{1}{\delta+\rho}\rho\alpha (Q+\rho I)^{-1}q.
\end{align*}
Consider the fixed-point candidate $z^\star$ satisfying $z^\star = E_R z^\star - \dfrac{1}{\delta+\rho}\rho\alpha(Q+\rho I)^{-1}q$ and $z^{k+1}-z^\star = E_R (z^{k}-z^\star)$.
The $z^k$-update in~\eqref{eqn:ADMM_L2_iterations_relaxation} converges (and so does the ADMM algorithm) if and only if the spectral radius of the error matrix in the above linear iterations is less than one. The eigenvalues of $E_R$ can be written as
\begin{align}
\zeta_R(\alpha,\rho,\lambda_i(Q))=
 1- \dfrac{\alpha \rho (\lambda_i(Q)+\delta)}{(\rho+\lambda_i(Q))(\rho+\delta)}.\label{eqn:zeta_expressions_relaxation}
\end{align}
Since $\rho, \delta,$ and $\lambda_i(Q) \in \R{}_{++}$, we see that
$0<\alpha<2\underset{i}{\min}\dfrac{(\rho+\delta)(\rho+\lambda_i(Q))}{\rho(\lambda_i(Q)+\delta)}$ implies that $\vert \zeta_R(\alpha, \rho, \lambda_i(Q))\vert<1$ for all $i$, which completes the first part of the proof.

For a fixed $\rho$ and $\delta$, we now characterize the values of $\alpha$ that ensure that the over-relaxed iterations~\eqref{eqn:ADMM_L2_iterations_relaxation} have a smaller convergence factor and thus a smaller $\varepsilon$-solution time than the classical ADMM iterates~\eqref{eqn:ADMM_L2_iterations}, \ie $\zeta_R-\zeta < 0$. From~\eqref{eqn:fexpressions} and~\eqref{eqn:zeta_expressions_relaxation} we have  $\operatorname{argmax}_i{\zeta_R(\alpha, \rho, \lambda_i(Q))} = \operatorname{argmax}_i{\zeta(\rho, \lambda_i(Q))}$, since $\zeta_R$ and $\zeta$ are equivalent up to an affine transformation and they have the same sign of the derivative with respect to $\lambda_i(Q)$.
For any given $\lambda_{i}(Q)$ we have
\begin{align*}
\zeta_R-\zeta = \dfrac{\rho(1-\alpha)(\lambda_i(Q)+\delta) }{\rho^2+(\lambda_i(Q)+\delta)\rho+\lambda_i(Q) \delta}
\end{align*}
and we conclude that $\zeta_R-\zeta <0$ when $\alpha\in \left(1,\;\dfrac{2(\rho+\delta)(\rho+\lambda_i(Q))}{\rho(\lambda_i(Q)+\delta)}\right)$. Recalling the first part of the proof we conclude that, for given $\rho,\delta\in\R{}_{++}$, the over-relaxed iterations converge with a smaller convergence factor than classical ADMM for $1<\alpha<2\underset{i}{\min}\dfrac{(\rho+\delta)(\rho+\lambda_i(Q))}{\rho(\lambda_i(Q)+\delta)}$.


To find $(\rho^\star, \alpha^\star, \zeta_R^\star)$, we define
\begin{equation}
\begin{aligned}
\label{eqn:ADMM_L2_optimal_step-size_optProblem_relaxation}
(\rho^{\star}, \alpha^\star) = \underset{\rho, \alpha}{\mbox{argmin }} \max_i\left\vert \zeta_R(\rho, \alpha,\lambda_i(Q))\right\vert, \quad
\zeta_R^{\star} = \max_i\left \vert \zeta_R(\rho^\star,\alpha^\star,\lambda_i(Q))\right\vert.
\end{aligned}
\end{equation}

One readily verifies that $\zeta_R(\delta, 2, \lambda_i(Q))=0$ for $i= 1,\dots n$. Since zero is the global minimum of $\vert \zeta_R\vert$ we conclude that the pair $(\rho^\star, \alpha^\star) = (\delta, 2)$ is optimal. Moreover, for $(\rho^\star, \alpha^\star) = (\delta, 2)$ the matrix $E_R$ is a matrix of zeros and thus the algorithm~\eqref{eqn:ADMM_L2_iterations_relaxation} converges in one iteration. 
\subsection{Proof of Proposition~\ref{prop:w_residuals}}
For the sake of brevity we derive the expressions only for $w_-^{k+1}\triangleq F^{k+1}v^{k+1}-F^kv^k$, as similar computations also apply to $w_+^{k+1}\triangleq v^{k+1}-v^k$.
First, since $v^k=z^k+u^k$, it holds that $F^k v^k = (2D^k - I)v^k = 2D^k v^k - u^k-z^k$. From the equality $D^kv^k = u^k$ we then have $F^kv^k = u^k-z^k$. The residual $w_-^{k+1} $ can be rewritten as $w_-^{k+1} = u^{k+1} - u^k - z^{k+1} + z^k$. From~\eqref{eq:primal_res} and~\eqref{eqn:Quadratic_admm_iterations} we observe that $u^{k+1}-u^k = r^{k+1}$, so $w_-^{k+1} = r^{k+1} - (z^{k+1} - z^k)$. Decomposing $z^{k+1} - z^k$ as $\Pi_{\Range{A}}(z^{k+1} - z^k) + \Pi_{\Null{A^\top}}(z^{k+1} - z^k)$ we then conclude that $w_-^{k+1} = r^{k+1} - \Pi_{\Range{A}}(z^{k+1} - z^k) - \Pi_{\Null{A^\top}}(z^{k+1} - z^k)$.
We now examine each case $(i)-(iii)$ separately:

(i) When $A$ has full column rank, $\Pi_{\Range{A}} = A(A^\top A)^{-1}A^\top$ and $\Pi_{\Null{A^\top}} = I-\Pi_{\Range{A}}$. In the light of the dual residual~\eqref{eq:dual_res} we obtain  $\Pi_{\Range{A}}(z^{k+1} - z^k) =1/\rho A(A^\top A)^{-1}s^{k+1}$.

(ii)  Note that the nullity of $A^\top$ is $0$ if $A$ is full row-rank. Thus, $\Pi_{\Null{A^\top}} = 0$ and $\Pi_{\Range{A}}=I$. Moreover, since $AA^\top$ is invertible, $z^{k+1}-z^{k} = (AA^\top)^{-1}AA^\top(z^{k+1}-z^{k}) = 1/\rho(AA^\top)^{-1}As^{k+1}$.

(iii) When $A$ is invertible, the result easily follows.

We now relate the norm of $ r^{k+1}$ and $ s^{k+1}$ to the one of $ w_-^{k+1}$.
From~\eqref{eqn:w_minus_residual} and~\eqref{eqn:w_plus_residual}, we have
\begin{align*}
  \Vert r^{k+1}\Vert = \dfrac{1}{2}\Vert w_-^{k+1} + w_+^{k+1} \Vert \leq \dfrac{1}{2} (\Vert w_-^{k+1} \Vert + \Vert w_+^{k+1}\Vert)\leq \Vert w_-^{k+1}\Vert,
\end{align*}
where the first inequality is the triangle inequality and the last inequality holds as $v^k$'s are positive vectors, $\Vert w_+^{k+1}\Vert =\Vert v^{k+1}-v^k\Vert\leq \Vert F^{k+1}v^{k+1}-F^k v^k \Vert = \Vert w_-^{k+1}\Vert$.

For the dual residual, it can be verified that in case (i) and (ii) $A^\top(w_+^{k+1} - w_-^{k+1}) = \dfrac{2}{\rho}s^{k+1}$, so
\begin{align*}
  \Vert s^{k+1}\Vert = &\dfrac{\rho}{2}\Vert A^\top (w_-^{k+1} - w_+^{k+1}) \Vert \leq \dfrac{\rho}{2} \Vert A\Vert\left(\Vert w_-^{k+1} -  w_+^{k+1}\Vert\right)\\
  &\leq \dfrac{\rho}{2} \Vert A\Vert\left(\Vert w_-^{k+1} \Vert + \Vert w_+^{k+1}\Vert\right)\leq \rho \Vert A\Vert \Vert w_-^{k+1}\Vert.
\end{align*}
In case (iii), one finds $A(w_+^{k+1} - w_-^{k+1})= \dfrac{2}{\rho}s^{k+1}$ and again the same bound can be achieved (by replacing $A^\top$ with $A$ in above equality), thus concluding the proof.

\subsection{Proof of Theorem~\ref{thm:QP:linear_rate}}
Note that since $v^k$ is positive and $F^{k}$ is diagonal with elements in $\pm 1$, $F^{k+1}v^{k+1}=F^{k}v^{k}$ implies $v^{k+1}=v^k$. Hence, it suffices to establish the convergence of $F^{k}v^k$. From~\eqref{eqn:QP_Fv_sequence} we have
\begin{align*}
 \left\Vert F^{k+1}v^{k+1} - F^{k} v^{k} \right\Vert \leq \dfrac{1}{2}\left\Vert 2M - I\right\Vert \left\Vert v^{k} - v^{k-1} \right\Vert + \dfrac{1}{2} \left\Vert F^k v^k - F^{k-1} v^{k-1}\right\Vert.
 \end{align*}
Furthermore, as  $v^k$s are positive vectors, $\left\Vert v^k - v^{k-1}\right\Vert\leq \left\Vert F^k v^k - F^{k-1} v^{k-1}\right\Vert$, which implies
\begin{align}
\label{eqn:Quadratic_linear_rate}
\left\Vert F^{k+1}v^{k+1} - F^k v^k \right\Vert \leq  \underset{\zeta}{\underbrace{\left( \dfrac{1}{2}\left\Vert 2M - I \right\Vert + \dfrac{1}{2} \right)}} \left\Vert F^k v^k - F^{k-1} v^{k-1} \right\Vert.
\end{align}
We conclude that if $\left\Vert 2M - I\right\Vert < 1$, then $\zeta<1$ and the iterations~\eqref{eqn:QP_Fv_sequence} converge to zero at a linear rate.

To determine for what values of $\rho$ the iterations~\eqref{eqn:QP_Fv_sequence} converge, we characterize the eigenvalues of $M$. By the matrix inversion lemma $M =   \rho AQ^{-1}A^\top - \rho AQ^{-1} A^\top (I+ \rho A Q^{-1} A^\top)^{-1} \rho A Q^{-1}A^\top$.
From~\cite[Cor. 2.4.4]{HoJ:85}, $(I+ \rho A Q^{-1} A^\top)^{-1}$ is a polynomial function of $\rho A Q^{-1} A^\top$ which implies that $M=f(\rho AQ^{-1}A^\top)$ is a polynomial function of $\rho AQ^{-1}A^\top$ with $f(t) = t - t (1+t)^{-1} t$. Applying~\cite[Thm. 1.1.6]{HoJ:85}, the eigenvalues of $M$ are given by $f(\lambda_i(\rho AQ^{-1}A^\top))$ and thus
 \begin{align}
 \label{eqn:QP_M_eigenvalues}
\lambda_i(M) = \dfrac{\lambda_i(\rho AQ^{-1}A^\top)}{1+ \lambda_i(\rho AQ^{-1}A^\top)}.
\end{align}

If $\rho > 0$, then $\lambda_i(\rho AQ^{-1}A^\top)\geq 0$ and $\lambda_i(M)\in [0,1)$. Hence $\left\Vert 2M - I\right\Vert  \leq 1$ is guaranteed for all $\rho \in \R{}_{++}$ and equality only occurs if $M$ has eigenvalues at $0$. If $A$ is invertible or has full row-rank, then $M$ is invertible and all its eigenvalues are strictly positive, so $\left\Vert 2M - I\right\Vert <1$ and~\eqref{eqn:QP_Fv_sequence} is guaranteed to converge linearly. The case when $A$ is tall, i.e., $A^\top$ is rank deficient, is more challenging since $M$ has zero eigenvalues and  $\left\Vert 2M - I\right\Vert  = 1$. To prove convergence in this case, we analyze the $0$-eigenspace of $M$ and show that it can be disregarded. From the $x$-iterates given in~\eqref{eqn:Quadratic_admm_iterations_reformed} we have ${x^{k+1}-x^{k} = -(Q/\rho+A^\top A)^{-1}A^\top (v^k - v^{k-1})}$. Multiplying the former equality by $A$ from the left on both sides yields ${A(x^{k+1}-x^{k}) = - M (v^k-v^{k-1})}$. Consider a nonzero vector ${v^k - v^{k-1}}$ in $\Null{M}$. Then we have either ${x^{k+1} = x^{k}}$ or ${x^{k+1}-x^{k}}\in \Null{A}$. Having assumed that $A$ is full column-rank denies the second hypothesis. In other words, the $0$-eigenspace of $M$ corresponds to the stationary points of the algorithm~\eqref{eqn:Quadratic_admm_iterations_reformed}. We therefore disregard this eigenspace and the convergence result holds.
Finally, the R-linear convergence of the primal and dual residuals follows from the linear convergence rate of $F^{k+1}v^{k+1}-F^kv^k$ and Proposition~\ref{prop:w_residuals}.
%
\subsection{Proof of Theorem~\ref{thm:QP_optimal_factor}}
From the proof of Theorem~\ref{thm:QP:linear_rate} recall that
 \begin{align*}
\left\Vert F^{k+1}v^{k+1} - F^k v^k \right\Vert \leq  \left( \dfrac{1}{2}\left\Vert 2M - I \right\Vert + \dfrac{1}{2} \right) \left\Vert F^k v^k - F^{k-1} v^{k-1} \right\Vert.
\end{align*}
Define
\begin{align*}
	\zeta \triangleq \frac{1}{2}\Vert 2M-I \Vert + \frac{1}{2} &= \max_i \frac{1}{2} \vert 2 \lambda_i(M)-1\vert + \frac{1}{2}
=
\max_i
\left\vert \dfrac{\rho \lambda_i(AQ^{-1}A^\top)}{1+\rho\lambda_i(A Q^{-1} A^\top)} - \dfrac{1}{2}\right\vert + \dfrac{1}{2}
\end{align*}
where the last equality follows from the definition of $\lambda_i(M)$ in~\eqref{eqn:QP_M_eigenvalues}. Since $\rho>0$ and for the case where $A$ is either invertible or has full row-rank, $\lambda_i(AQ^{-1}A^{\top}) >0$ for all $i$, we conclude that $\zeta<1$.

It remains to find $\rho^\star$ that minimizes the convergence factor, \emph{i.e.}
\begin{align}\label{eqn:QP_optProblem}
\rho^{\star}
&= \underset{\rho}{\mbox{argmin }} \max_i\left\{ \left\vert \dfrac{\rho \lambda_i(AQ^{-1}A^\top)}{1+\rho\lambda_i(A Q^{-1} A^\top)} - \dfrac{1}{2}\right\vert + \dfrac{1}{2}\right\}.
\end{align}
Since $\dfrac{\rho\lambda_i(A Q^{-1} A^\top)}{1+\rho\lambda_i(A Q^{-1} A^\top)}$ is a monotonically increasing function in $\lambda_i(A Q^{-1} A^\top)$, the maximum values of $\zeta$ happen for the two extreme eigenvalues $\lambda_1(A Q^{-1} A^\top)$ and $\lambda_n(A Q^{-1} A^\top)$:
%
\begin{equation}
\begin{aligned}\label{eqn:QP_convergence_factor}
\max_i \left\{ \zeta(\lambda_i(A Q^{-1} A^\top),\rho)\right\} =
\left\{
\begin{array}[c]{lll}
\dfrac{1}{1+\rho\lambda_1(A Q^{-1} A^\top)} & \mbox{if} &\rho \leq \rho^\star, \\
\dfrac{\rho\lambda_n(A Q^{-1} A^\top)}{1+\rho\lambda_n(A Q^{-1} A^\top)} & \mbox{if}& \rho > \rho^\star. \end{array}\right.
\end{aligned}
\end{equation}
Since the left brace of $\max_i \left\{ \zeta(\lambda_i(A Q^{-1} A^\top),\rho)\right\}$, i.e.~$\dfrac{1}{1+\rho\lambda_1(A Q^{-1} A^\top)}$ is monotone decreasing in $\rho$ and the right brace is monotone increasing, the minimum with respect to $\rho$ happens at the intersection point~\eqref{eqn:QP_optimal_factor}.
\subsection{Proof of Theorem~\ref{thm:QP_slow}}
First we derive the lower bound on the convergence factor and show it is strictly smaller than $1$. From~\eqref{eqn:QP_Fv_sequence} we have $\left\Vert F^{k+1}v^{k+1} - F^k v^k \right\Vert
= \left\Vert D^k v^k-D^{k-1} v^{k-1}   -  M (v^k - v^{k-1})\right\Vert.$
By applying the reverse triangle inequality and dividing by $\|F^kv^k - F^{k-1} v^{k-1}\|$, we find
\begin{equation*}
\dfrac{\|F^{k+1}v^{k+1} - F^k v^k \|}{\|F^k v^k - F^{k-1} v^{k-1}\|} \geq \vert \delta_k - \epsilon_k \vert.
\end{equation*}
Recalling from~\eqref{eqn:convergence_factor_def} that the convergence factor $\zeta$ is the maximum over $k$ of the left hand-side yields the lower bound~\eqref{eq:QP_lower_bound}. Moreover, the inequality $1 > \zeta \geq \underline{\zeta}$ follows directly from Theorem~\ref{thm:QP:linear_rate}.

The second part of the proof addresses the cases (i)-(iii) for $\rho>0$. Consider case (i) and let
$\Null{A^\top}= \{0\}$. It follows from Theorem~\ref{thm:QP_optimal_factor} that the convergence factor is given by $\tilde{\zeta}(\rho)$, thus proving the sufficiency of $\Null{A^\top}= \{0\}$ in (i). The necessity follows directly from statement (iii), which is proved later.

Now consider the statement (ii) and suppose $\Null{A^\top}$ is not zero-dimensional. Recall that $\lambda_1(AQ^{-1}A^\top)$ is the smallest nonzero eigenvalue of $AQ^{-1}A^\top$ and suppose that $\epsilon_k \geq 1-\xi$. Next we show that $\epsilon_k \geq 1-\xi$ implies $\|\Pi_{\Null{A^\top}}(v^k - v^{k-1})\| / \|v^k - v^{k-1}\|\leq \sqrt{2\xi}$. Since $M\Pi_{\Null{A^\top}} = 0$, $\|M\|<1$, and $\|v^k - v^{k-1}\|\leq \|F^kv^k - F^{k-1}v^{k-1}\|$ we have
\begin{equation*}
\begin{aligned}
\epsilon_k^2 &= \dfrac{\|M(I-\Pi_{\Null{A^\top}})( v^k - v^{k-1})  \|^2}{\|F^kv^k - F^{k-1}v^{k-1}\|^2}
\leq \dfrac{\|\Pi_{\Range{A}}(v^k - v^{k-1})  \|^2}{\|v^k - v^{k-1}\|^2}
= 1-\dfrac{\|\Pi_{\Null{A^\top}}(v^k - v^{k-1})  \|^2}{\|v^k - v^{k-1}\|^2}.
\end{aligned}
\end{equation*}
Using the above inequality and $\epsilon_k^2 \geq (1-\xi)^2$ we obtain $\|\Pi_{\Null{A^\top}}(v^k - v^{k-1})\| / \|v^k - v^{k-1}\|\leq \sqrt{2\xi-\xi^2} \leq \sqrt{2\xi}$.

The latter inequality allows us to derive an upper-bound on $\underline{\zeta}$ as follows.
Recalling~\eqref{eqn:QP_Fv_sequence}, we have
\begin{equation}\label{eq:QP_slow_ii}
\begin{aligned}
\underline{\zeta}\leq\dfrac{\|F^{k+1}v^{k+1} - F^k v^k \|}{\|F^k v^k - F^{k-1} v^{k-1}\|}
&\leq
 \dfrac{1}{2} + \dfrac{1}{2}\dfrac{\|(I- 2M) (v^k - v^{k-1})\|}{\|F^k v^k - F^{k-1} v^{k-1}\|} \\
&= \dfrac{1}{2} + \dfrac{1}{2}\sqrt{\dfrac{\|(I- 2M)\Pi_{\Range{A}}(v^k - v^{k-1})\|^2}{\|F^k v^k - F^{k-1} v^{k-1}\|^2} + \dfrac{\|\Pi_{\Null{A^\top}}(v^k - v^{k-1})\|^2}{\|F^k v^k - F^{k-1} v^{k-1}\|^2}}.
 \end{aligned}
\end{equation}
Using the inequalities $\|v^k - v^{k-1}\| \leq \|F^k v^k - F^{k-1} v^{k-1}\|$ and $\sqrt{a^2 + b^2} \leq a + b$ for $a,b\in\mathcal{R}_+$, the inequality~\eqref{eq:QP_slow_ii} becomes $\underline{\zeta}\leq \dfrac{1}{2} + \dfrac{1}{2}\|(I- 2M)\Pi_{\Range{A}}\| + \sqrt{\dfrac{\xi}{2}} \leq  \tilde{\zeta}(\rho) + \sqrt{\dfrac{\xi}{2}},$
which concludes the proof of (ii).

As for the third case (iii), note that $\epsilon_k\leq \xi$ holds if $\|\Pi_{\Null{A^\top}}(v^k - v^{k-1})\| / \|v^k - v^{k-1}\|\geq \sqrt{1-\xi^2/\|M\|^2}$, as the latter inequality implies that
\begin{equation*}
\begin{aligned}
\epsilon_k &= \dfrac{\|M\Pi_{\Range{A}}(v^k - v^{k-1})\|}{\|F^{k}v^k - F^{k-1}v^{k-1}\|} \leq \|M\|\dfrac{\|\Pi_{\Range{A}}(v^k - v^{k-1})\|}{\|v^k - v^{k-1}\|} \leq \xi.
\end{aligned}
\end{equation*}
Supposing that there exists a non-empty set $\mathcal{K}$ such that $\delta_k \geq 1-\xi$ and
$\|\Pi_{\Null{A^\top}}(v^k - v^{k-1})\| / \|v^k - v^{k-1}\|\geq \sqrt{1-\xi^2/\|M\|^2}$ holds for all $k\in\mathcal{K}$, we have
$\underline{\zeta} \geq \max_{k\in\mathcal{K}}\; \delta_k - \epsilon_k \geq 1-2\xi$ regardless the choice of $\rho$.
\subsection{Proof of Lemma~\ref{lem:optimal_fixed_point}}
Let $(x^\star,\, z^\star,\, u^\star)$ denote a fixed-point of~\eqref{eqn:Quadratic_admm_iterations_relaxation} and let $\mu$ be the Lagrange
multiplier associated with the equality constraint in~\eqref{eqn:Quadratic_problem_1}. For the optimization problem~\eqref{eqn:Quadratic_problem_1}, the Karush-Kuhn-Tucker (KKT) optimality conditions~\cite{Nesterov03} are
\begin{equation*}
\begin{aligned}
0 & = Q x + q + A^\top \mu,\quad\;\; z \geq 0,\\
0 &=Ax + z - b,\quad\quad\quad 0  =\mbox{diag}(\mu) z.
\end{aligned}
\end{equation*}
Next we show that the KKT conditions hold for the fixed-point $(x^\star,\, z^\star,\, u^\star)$ with $\mu^\star = 1/\rho u^\star$.
From the $u-$iterations we have $0 = \alpha (A x^\star-c)-(1-\alpha)z^\star+z^\star = \alpha( A x^\star +z^\star-c)$. It follows that $z^\star$ is given by $z^{\star} = \mbox{max}\{0,-\alpha (A x^{\star} + z^\star -c)+z^\star-u^{\star}\} = \mbox{max}\{0, z^\star-u^{\star}\} \geq 0$. The $x-$iteration then yields
%
$0 = Qx^\star +q + \rho A^\top(Ax^\star + z^\star -c + u^\star) = Qx^\star +q + A^\top \mu^\star$.
Finally,  from
$z^\star \geq 0$ and the $z-$update, we have that $z^\star_i > 0 \Rightarrow u^\star_i = 0$ and $z^\star_i = 0 \Rightarrow u^\star_i \geq 0$. Thus, $\rho\,\mbox{diag}(\mu^\star) z^\star = 0$.
%
%
%

\subsection{Proof of Theorem~\ref{thm:QP_relaxation_convergence}}
 Taking the Euclidean norm of~\eqref{eqn:QP_Fv_sequence_relaxation} and applying the Cauchy-Schwarz inequality yields
 \begin{align*}
 \left\Vert F^{k+1}v^{k+1} - F^{k} v^{k} \right\Vert \leq  \dfrac{\vert \alpha \vert}{2}\left\Vert 2M - I \right\Vert \left\Vert v^{k} - v^{k-1} \right\Vert + \vert 1-\dfrac{\alpha}{2}\vert \left\Vert F^k v^k - F^{k-1} v^{k-1}\right\Vert.
 \end{align*}
Note that since $v^k$s are positive vectors we have $\left\Vert v^k - v^{k-1}\right\Vert\leq \left\Vert F^k v^k - F^{k-1} v^{k-1}\right\Vert$ and thus
\begin{align}
\label{eqn:Quadratic_linear_rate_relaxation}
\dfrac{\left\Vert F^{k+1}v^{k+1} - F^{k} v^{k} \right\Vert}{\left\Vert F^k v^k - F^{k-1} v^{k-1} \right\Vert} \leq  \underset{\zeta_R}{\underbrace{\left(  \dfrac{\vert \alpha \vert}{2}\left\Vert 2M - I \right\Vert + \left\vert 1-\dfrac{\alpha}{2} \right\vert \right)}} .
\end{align}

Note that $\rho \in \R{}_{++}$ and recall from the proof of Theorem~\ref{thm:QP:linear_rate} that the $0$-eigenspace of $M$ can be disregarded. Therefore, $\dfrac{1}{2}\left\Vert 2M-I\right\Vert_{\Null{M}^\bot} \in [0,\dfrac{1}{2})$. 
%
%
Defining $\tau \triangleq \dfrac{1}{2}\left\Vert 2M- I\right\Vert_{\Null{M}^\bot}$ 
we have \[ \zeta_R=\alpha \tau + \vert1- \dfrac{\alpha}{2}\vert <  \dfrac{\alpha}{2} + \vert1- \dfrac{\alpha}{2}\vert \]
Hence, we conclude that for $\rho \in \R{}_{++}$ and $\alpha\in (0, 2]$, it holds that $\zeta_R <1$ , which implies that~\eqref{eqn:QP_Fv_sequence_relaxation} converges linearly to a fixed-point. By  Lemma~\ref{lem:optimal_fixed_point} this fixed-point is also a global optimum of \eqref{eqn:Quadratic_problem}.  Now, denote $w_-^{k+1} \triangleq F^{k+1} v^{k+1}- F^{k} v^{k}$ and $w_+^{k+1}\triangleq v^{k+1}-v^k$. Following the same steps as Proposition~\ref{prop:w_residuals}, it is easily verified that $w_-^{k+1} = u^{k+1}-u^k + z^k - z^{k+1}$ and $w_+^{k+1}=u^{k+1}-u^k +z^{k+1}- z^k$ from which combined with~\eqref{eqn:Quadratic_admm_iterations_relaxation} one obtains
\begin{align*}
 s^{k+1} = \rho\dfrac{A^\top}{2} (w_+^{k+1}-w_-^{k+1}), \quad r^{k+1} = \dfrac{1}{2}w_+^{k+1} + \dfrac{2-\alpha}{2\alpha}w_-^{k+1}.
\end{align*}
We only upper-bound $\Vert r^{k+1}\Vert$, since an upper bound for $\Vert s^{k+1}\Vert$ was already established in~\eqref{eqn:s_to_fv}. Taking the Euclidean norm of the second equality above and using the triangle inequality
\begin{align}\label{eqn:r_to_fv_relaxed}
\Vert r^{k+1}\Vert \leq  \dfrac{1}{2}  \Vert w_+^{k+1}\Vert + \dfrac{2-\alpha}{2\alpha} \Vert w_-^{k+1}\Vert \leq  
\dfrac{1}{\alpha} \Vert w_-^{k+1}\Vert.
\end{align}
The R-linear convergence of the primal and dual residuals now follows from the linear convergence rate of $F^{k+1}v^{k+1}-F^kv^k$ and the bounds in~\eqref{eqn:s_to_fv} and~\eqref{eqn:r_to_fv_relaxed}.
%


\subsection{Proof of Theorem~\ref{thm:QP_relaxation_optimal_factor}}
Define
\begin{equation}
\begin{aligned}
\label{eqn:QP_Convergencefactor_relaxation}
&\zeta_R( \rho, \alpha ,\lambda_i(AQ^{-1}A^\top)) = \alpha \left\vert \dfrac{\rho \lambda_i(AQ^{-1}A^\top)}{1 + \rho \lambda_i(AQ^{-1}A^\top)} - \dfrac{1}{2}\right\vert + 1- \dfrac{\alpha}{2},\\
&\zeta_R^\star = \underset{i}{\max} \, \underset{\rho, \alpha}{\min}\{\zeta_R(\rho, \alpha ,\lambda_i(AQ^{-1}A^\top))\}.
\end{aligned}
\end{equation}
Since $\left\vert \dfrac{\rho \lambda_i(AQ^{-1}A^\top)}{1 + \rho \lambda_i(AQ^{-1}A^\top)} - \dfrac{1}{2}\right\vert < \dfrac{1}{2}$, it follows that $\zeta_R(\rho, \alpha, \lambda_i(AQ^{-1}A^\top))$ is monotone decreasing in $\alpha$. Thus, $\zeta_R(\rho, \alpha, \lambda_i(AQ^{-1}A^\top))$ is minimized by $\alpha^\star=2$. To determine
\begin{align}
\label{eqn:QP_optProblem_relaxation}
\rho^{\star} = \underset{\rho}{\mbox{argmin }} \max_i\left\{ \zeta_R(\rho, 2,\lambda_i(AQ^{-1}A^\top))\right\},
\end{align}

we note that~\eqref{eqn:QP_optProblem} and~\eqref{eqn:QP_optProblem_relaxation} are equivalent up to an affine transformation, hence we have the same minimizer $\rho^\star$. It follows from the proof of Theorem~\ref{thm:QP_optimal_factor} that ${\rho^\star = 1/\sqrt{\lambda_1(AQ^{-1}A^\top)\; \lambda_n(AQ^{-1}A^\top)}}$. Using $\rho^\star$ in~\eqref{eqn:QP_Convergencefactor_relaxation} results in the convergence factor~\eqref{eqn:QP_relaxation_optimal_factor}.

For given $A$, $Q$, and $\rho$, we can now find the range of values of $\alpha$ for which~\eqref{eqn:Quadratic_admm_iterations_relaxation} have a smaller convergence factor than~\eqref{eqn:Quadratic_admm_iterations}, \ie for which $\zeta_R - \zeta < 0$. By~\eqref{eqn:Quadratic_linear_rate} and~\eqref{eqn:Quadratic_linear_rate_relaxation} it holds that
\begin{align*}
\zeta_R - \zeta =  \dfrac{\alpha}{2}\left\Vert 2M - I \right\Vert + 1-\dfrac{\alpha}{2}- \dfrac{1}{2}\left\Vert 2M - I\right\Vert - \dfrac{1}{2} =  (1 - \alpha )\left(\dfrac{1}{2} - \dfrac{1}{2}\left\Vert 2M - I\right\Vert  \right).
\end{align*}
This means that $\zeta_R - \zeta<0$ when $\alpha>1$. Therefore, the iterates produced by the relaxed algorithm~\eqref{eqn:Quadratic_admm_iterations_relaxation} have smaller convergence factor than the iterates produced by~\eqref{eqn:Quadratic_admm_iterations}
for all values of the relaxation parameter $\alpha \in (1,2]$. This concludes the proof.
\subsection{Proof of Theorem~\ref{thm:QP_optimal_preconditioning}}
Note that the non-zero eigenvalues of $LAQ^{-1}A^\top L$ are the same as the ones of $R_q^\top A^\top W A R_q$ where $W=L^2$ and $R_q^\top R_q=Q^{-1}$ is its Choleski factorization~\cite{HoJ:85}. Defining $\lambda_n(R_q^\top A^\top W A R_q)$ and $\lambda_1(R_q^\top A^\top W A R_q)$ as the largest and smallest nonzero eigenvalues of $LAQ^{-1}A^\top L$, the optimization problem we aim at solving can be formulated as
\begin{equation}\label{eqn:QP_optimal_scaling_proof}
\begin{aligned}
\begin{array}{ll}
\underset{\bar\lambda\in\R{},\;\underline\lambda\in\R{},\;l\in\R{m}}{\mbox{minimize}} & {\bar\lambda} / {\underline\lambda}\\
\mbox{subject to}& \bar\lambda > \lambda_n(R_q^\top A^\top W A R_q),\\
& \lambda_1(R_q^\top A^\top W A R_q) > \underline\lambda,\\
& W=\mbox{diag}(w),\; w > 0.
\end{array}
\end{aligned}
\end{equation}
In the proof we show that the optimization problem \eqref{eqn:QP_optimal_scaling_proof} is equivalent to~\eqref{eqn:QP_optimal_scaling_convex}.

Define $T(\bar\lambda) \triangleq \bar\lambda I - R_q^\top A^\top W A R_q$. First observe that $\bar\lambda \geq \lambda_n(R_q^\top A^\top W A R_q)$ holds if and only if $T(\bar\lambda) \in \PSD{n}$, which proves the first inequality in the constraint set~\eqref{eqn:QP_optimal_scaling_convex}.


To obtain a lower bound on $\lambda_1(R_q^\top A^\top W A R_q)$ one must disregard the zero eigenvalues of $R_q^\top A^\top W A R_q$ (if they exist). This can be performed by restricting ourselves to  the subspace orthogonal to $\Null{R_q^\top A^\top W A R_q} = \Null{A R_q}$. In fact, letting $s$ to be the dimension of the nullity of $A R_q$ or simply $A$ and denoting $P^{n\times n-s}$ as a basis of $\mbox{Im}(R_q^\top A^\top)$, we have that $\underline\lambda \leq \lambda_1$ if and only if $x^\top P^\top  T(\underline\lambda) Px \leq 0$ for all $x\in\R{n-s}$. Note that for the case when the nullity of $A$ is $0$ ($s=0$), all the eigenvalues of $R_q^\top A^\top W A R_q$ are strictly positive and, hence, one can set $P=I$. We conclude that $\underline\lambda \leq \lambda_1(R_q^\top A^\top W A R_q)$ if and only if $P^\top\left(R_q^\top A^\top W A R_q-\underline\lambda I \right) P \in \PSD{n-s}$.

Note that $\lambda_1(R_q^\top A^\top W A R_q)>0$ can be chosen arbitrarily by scaling $W$, which does not affect the ratio $\lambda_n(R_q^\top A^\top W A R_q) / \lambda_1(R_q^\top A^\top W A R_q)$. Without loss of generality, one can suppose $\underline\lambda^\star = 1$ and thus the lower bound on $\lambda_1(R_q^\top A^\top W A R_q)\geq \underline\lambda^\star = 1$ corresponds to the last inequality in the constraint set of~\eqref{eqn:QP_optimal_scaling_convex}. Observe that the optimization problem now reduces to minimizing $\bar\lambda$.
The proof concludes by rewriting~\eqref{eqn:QP_optimal_scaling_proof} as~\eqref{eqn:QP_optimal_scaling_convex}, which is a convex problem.
\subsection{Proof of Proposition~\ref{prop:QP_when_careful_1}}
Assuming ${F^{k+1} = F^{k} = -I}$,~\eqref{eqn:v_recurrence_relaxation} reduces to
$v^{k+1}-v^k = \left( (1-\alpha)I + \alpha M\right) (v^k-v^{k-1})$.
By taking the Euclidean norm of both sides and applying the Cauchy inequality, we find
 \begin{align*}
    \Vert v^{k+1} - v^{k} \Vert \leq \Vert (1-\alpha)I + \alpha M \Vert \Vert v^k-v^{k}\Vert.
 \end{align*}
Since the eigenvalues $M$ are $\dfrac{\rho \lambda_i(AQ^{-1}A^\top)}{1+\rho \lambda_i(AQ^{-1}A^\top)}$, the convergence factor $\zeta_R$ is
\begin{align*}
    \zeta_R(\rho,\alpha,\lambda_i(AQ^{-1}A^\top)) &= 1-\alpha + \alpha \dfrac{\rho \lambda_i(AQ^{-1}A^\top)}{1+\rho \lambda_i(AQ^{-1}A^\top)}. 
\end{align*}
It is easy to check that the smallest value of $\vert \zeta_R\vert$ is obtained when $\alpha=1$ and $\rho\rightarrow 0$. Since $\alpha=1$ the relaxed ADMM iterations~\eqref{eqn:Quadratic_admm_iterations_relaxation} coincide with~\eqref{eqn:Quadratic_admm_iterations} and consequently $\zeta=\zeta_R$.
\subsection{Proof of Proposition~\ref{prop:QP_when_careful_2}}
The proof follows similarly to the one of Proposition~\ref{prop:QP_when_careful_1} but with ${F^{k+1} = F^{k} = I}$.

\end{document}